\documentclass[11pt]{article}

\IfFileExists{lmodern.sty}{%
  \usepackage[T1]{fontenc}
  \usepackage{lmodern}
}{}
\usepackage{amsmath}
\usepackage{amssymb}
\usepackage{bm}
\usepackage{booktabs}
\usepackage{tabularx}
\usepackage{graphicx}
\usepackage{microtype}
\usepackage{siunitx}
\usepackage{hyperref}

\hypersetup{
  pdftitle={Hierarchical Log-Gaussian Relaxation on a Fixed D3Q125 Velocity Set},
  pdfauthor={Bjørn Wu},
  pdfsubject={Order-resolved hierarchical relaxation for a fixed D3Q125 discrete-velocity kinetic model},
  pdfkeywords={discrete velocity method, D3Q125, Hermite expansion, thermodynamic nonequilibrium, hierarchical relaxation, finite-volume kinetic method}
}

\title{Hierarchical Log-Gaussian Relaxation on a Fixed D3Q125 Velocity Set}
\author{%
Bjørn Wu\\
\small Independent Researcher, Oslo, Norway
}
\date{}

\begin{document}

\maketitle

\begin{abstract}
A hierarchical order-resolved relaxation model is developed for a fixed D3Q125 discrete-velocity kinetic formulation. Conventional adaptive collision models often use one scalar rarefaction or nonequilibrium indicator for all retained moment orders, thereby coupling otherwise distinct kinetic sectors. The present model combines a shared macroscopic-gradient background with separate second-, third-, and fourth-order thermodynamic nonequilibrium indicators. The resulting effective measures \(K_2\), \(K_3\), and \(K_4\) activate distinct log-Gaussian relaxation spectra without changing the velocity set.

Pure-order perturbation tests verify selective activation with nonmatching moment sectors remaining at roundoff level. Homogeneous mixed-order, amplitude, and composition tests show that the order-resolved formulation reduces residual nonequilibrium relative to a common-sensor model while maintaining positive populations. In a smooth periodic compression wave at the stated reference discretization and in the TNE-only sensor limit, the peak total nonequilibrium intensity is reduced by \(6.565\%\). The reduction occurs throughout the domain and in each retained moment sector, with the largest relaxation correction in the fourth-order channel.

Additional fixed-grid timestep, transport-discretization, relaxation-spectrum, uniform-boost, long-time, and shear-wave studies show that the sign of the hierarchical correction is robust over the tested configurations, whereas its magnitude remains dependent on timestep, transport scheme, sensor frame, and prescribed relaxation spectrum. Quantitative boost checks identify the laboratory-frame raw-Hermite origin of the frame sensitivity without detecting an evident collision-path inconsistency over the tested range. The periodic benchmarks preserve the principal global invariants to floating-point accuracy and remain positive. The results establish the mechanism, selectivity, and numerical behavior of order-resolved activation on a fixed velocity set; independent kinetic-reference validation is still required before claiming universal accuracy improvement.
\end{abstract}

\section{Introduction}

Kinetic descriptions provide access to flow information that is hidden when the particle distribution is reduced directly to a small set of macroscopic conservation laws. This additional information becomes important whenever the local distribution departs appreciably from equilibrium, including rarefied transport, thermal nonequilibrium, high-speed compressible flow, and under-resolved multiscale dynamics.

Discrete-velocity and lattice kinetic methods approximate the velocity dependence of the distribution by a finite set of populations. A systematic connection between discrete velocity sets, Hermite expansions, and Gauss--Hermite quadrature makes it possible to construct models with controlled moment accuracy. In this framework, increasing the quadrature order and retaining higher Hermite coefficients extends the range of kinetic moments represented by the discrete model \cite{grad1949kinetic,shan1998discretization,shan2006kinetic,shan2010general}.

The presence of higher-order moments creates a corresponding modeling question for the collision operator. A single-relaxation-time model forces all nonconserved components toward equilibrium at the same rate. Multiple-relaxation-time and central-moment formulations address this limitation by assigning distinct relaxation parameters to different moment channels. Such formulations can improve stability, separate transport processes, and reduce undesirable coupling among physical and nonhydrodynamic modes \cite{dhumieres2002multiple,lallemand2000theory,lallemand2003thermal,shan2019centralmoment,premnath2012centralmoment}.

High-order regularized lattice Boltzmann methods provide another important perspective. By reconstructing the nonequilibrium distribution from selected Hermite coefficients, regularization can filter unresolved kinetic content and substantially improve numerical stability \cite{coreixas2017recursive}. However, stability analyses also show that equilibration or filtering of high-order moments can modify the dissipation of modes carrying physical waves \cite{masset2020hydrodynamics,wissocq2020stability}. Consequently, higher-order moments should not be regarded merely as numerical noise: their treatment determines both stability and the kinetic response of the model.

Discrete Boltzmann modeling further emphasizes the diagnostic role of nonequilibrium moments. Deviations of stress-, heat-flux-, and higher-order moments from their equilibrium values can be used to characterize thermodynamic nonequilibrium and distinguish different stages or structures in multiscale flows \cite{gan2015nonequilibrium,guo2025threeDimensional}. These diagnostics retain information that is lost when the complete nonequilibrium state is compressed into a single scalar measure.

Despite these developments, adaptive relaxation models commonly face a remaining aggregation problem. A scalar rarefaction or total nonequilibrium indicator is often evaluated first and then supplied to all retained moment orders. In such a construction, a strongly nonequilibrium contribution in one sector changes the relaxation of every other sector. The resulting coupling is not imposed by the moment representation itself; it is introduced by the scalar sensor.

This issue is especially relevant when the retained moment orders represent qualitatively different departures from equilibrium. A second-order deviation is associated primarily with stress-like nonequilibrium, a third-order deviation contains heat-flux-like information, and fourth-order moments describe still higher kinetic structure. Their amplitudes need not be proportional, and they need not become important at the same spatial location or physical scale.

The present work introduces a hierarchical order-resolved activation model that retains this distinction. A common macroscopic rarefaction background is constructed from normalized density, temperature, and velocity gradients. It is then combined separately with the second-, third-, and fourth-order thermodynamic nonequilibrium contributions to form three effective indicators,

\[
K_2,\qquad K_3,\qquad K_4.
\]

Each indicator drives its own order-dependent log-Gaussian relaxation curve. This schedule adapts the log-Knudsen scale-space probability partition previously introduced for continuum--ballistic flux blending \cite{wu2026scaleSpace} to moment-specific collision factors. The model therefore preserves common activation by macroscopic gradients while avoiding premature aggregation of the kinetic nonequilibrium channels.

The velocity representation is held fixed at D3Q125 throughout this study. This choice isolates the effect of hierarchical collision activation from the separate problem of adaptive velocity-space refinement. No switching of velocity sets, local velocity scaling, or population remapping is used.

The principal contributions of this work are as follows:

\begin{enumerate}
\item A retained-subspace distribution-level decomposition into equilibrium and second-, third-, and fourth-order nonequilibrium Hermite sectors on a fixed D3Q125 velocity set.
\item A self-consistent rarefaction sensor combining a common macroscopic-gradient background with order-specific thermodynamic nonequilibrium measures.
\item Three hierarchical effective Knudsen indicators that activate distinct order-dependent log-Gaussian relaxation spectra.
\item Positivity, conservation, entropy, and transport diagnostics integrated into both local and spatially resolved collision updates.
\item Numerical verification using pure-order perturbations, homogeneous mixed-order relaxation, amplitude and composition scans, a smooth temperature wave, and a multi-gradient smooth compression wave.
\item Direct common-versus-order-resolved comparisons demonstrating both a TNE-sensitive regime with measurable hierarchical corrections and a macroscopic-gradient-dominated regime in which the two formulations converge.
\end{enumerate}

The numerical results show that the model is selectively activated, remains positive in all principal benchmarks, and reduces artificial cross-order suppression of relaxation. In the smooth compression-wave test at zero macroscopic-gradient sensor scale (the TNE-only sensor limit), the order-resolved formulation lowers the peak total thermodynamic nonequilibrium intensity by a measurable amount, with reductions occurring throughout the domain and in every retained moment sector. As the shared macroscopic-gradient background increases, the common and order-resolved models approach one another smoothly.

The present study is intended to establish the mechanism and numerical behavior of order-resolved activation. A reduction in thermodynamic nonequilibrium is not by itself a proof of improved physical accuracy. Studies comparing high-order lattice kinetic models with independent discrete-velocity solutions show that accuracy depends on both the velocity-space representation and the spatial-temporal discretization \cite{shi2021accuracy}. External validation against an independent kinetic reference solution is therefore identified as the principal remaining step before a universal accuracy claim can be made.

The remainder of the article is organized as follows. Section 2 presents the D3Q125 discrete kinetic and Hermite formulation. Section 3 defines the hierarchical collision model and order-resolved rarefaction indicators. Section 4 describes the finite-volume transport scheme, time integration, periodic treatment, and numerical diagnostics. Section 5 reports the mechanism and spatial benchmarks. Section 6 discusses the physical interpretation, regime dependence, numerical robustness, and limitations. Section 7 summarizes the conclusions.

\section{Discrete kinetic formulation}

\subsection{Discrete-velocity representation}

The kinetic state is represented by a finite set of populations

\[
f_i(\boldsymbol{x},t),
\qquad
i=1,\ldots,125,
\]

associated with the three-dimensional discrete velocities

\[
\boldsymbol{\xi}_i
=
\left(
\xi_{i,x},
\xi_{i,y},
\xi_{i,z}
\right).
\]

The velocity set is constructed as the tensor product of a five-node one-dimensional Gauss--Hermite quadrature,

\[
\mathcal{V}_{125}
=
\mathcal{V}_{5}
\otimes
\mathcal{V}_{5}
\otimes
\mathcal{V}_{5},
\]

and therefore contains

\[
5^3=125
\]

discrete velocities.

The one-dimensional nodes and weights correspond to the normalized standard-Gaussian weight

\[
\omega_1(\xi)
=
\frac{1}{\sqrt{2\pi}}
\exp\!\left(-\frac{\xi^2}{2}\right).
\]

The associated tensor-product weight is

\[
\omega(\boldsymbol{\xi})
=
\frac{1}{(2\pi)^{3/2}}
\exp\!\left(-\frac{|\boldsymbol{\xi}|^2}{2}\right).
\]

The five-node rule is exact for one-dimensional polynomials through degree nine. The tensor-product rule inherits this exactness independently in each Cartesian coordinate.

The corresponding three-dimensional quadrature weights are products of the one-dimensional weights,

\[
w_i
=
w_{i_x}^{(1)}
w_{i_y}^{(1)}
w_{i_z}^{(1)}.
\]

The numerical abscissae and weights are generated directly by the implementation from the underlying one-dimensional quadrature, thereby avoiding discrepancies between the manuscript and the executable model.

The discrete kinetic equation is written as

\[
\frac{\partial f_i}{\partial t}
+
\boldsymbol{\xi}_i
\cdot
\nabla f_i
=
\Omega_i,
\]

where \(\Omega_i\) denotes the hierarchical collision operator introduced in Section 3.

\subsection{Macroscopic fields}

The local density is evaluated from the zeroth discrete moment,

\[
\rho
=
\sum_i f_i.
\]

The momentum density and flow velocity are

\[
\rho\boldsymbol{u}
=
\sum_i
\boldsymbol{\xi}_i f_i,
\]

and

\[
\boldsymbol{u}
=
\frac{
\sum_i\boldsymbol{\xi}_i f_i
}{
\rho
}.
\]

Let

\[
\boldsymbol{c}_i
=
\boldsymbol{\xi}_i-\boldsymbol{u}
\]

denote the peculiar velocity. In three velocity dimensions, the local temperature is obtained from the trace of the second central moment,

\[
T
=
\frac{1}{3\rho}
\sum_i
\left|
\boldsymbol{c}_i
\right|^2
f_i,
\]

in the nondimensional convention used by the present implementation.

The pressure associated with the local equilibrium state is therefore

\[
p=\rho T.
\]

All numerical benchmarks use positive density and temperature fields, and the collision diagnostics monitor whether the reconstructed populations remain positive.

\subsection{Central moments}

The rank-\(n\) central moment tensor is defined as

\[
\boldsymbol{M}^{(n)}
=
\sum_i
f_i
\boldsymbol{c}_i^{\otimes n},
\]

where

\[
\boldsymbol{c}_i^{\otimes n}
=
\underbrace{
\boldsymbol{c}_i
\otimes\cdots\otimes
\boldsymbol{c}_i
}_{n\ \mathrm{times}}.
\]

The second central moment is

\[
M_{\alpha\beta}^{(2)}
=
\sum_i
f_i
c_{i,\alpha}
c_{i,\beta},
\]

the third central moment is

\[
M_{\alpha\beta\gamma}^{(3)}
=
\sum_i
f_i
c_{i,\alpha}
c_{i,\beta}
c_{i,\gamma},
\]

and the fourth central moment is

\[
M_{\alpha\beta\gamma\delta}^{(4)}
=
\sum_i
f_i
c_{i,\alpha}
c_{i,\beta}
c_{i,\gamma}
c_{i,\delta}.
\]

At local equilibrium, the odd central moments vanish and the second central moment is isotropic,

\[
M_{\alpha\beta,\mathrm{eq}}^{(2)}
=
\rho T\delta_{\alpha\beta}.
\]

The equilibrium fourth central moment has the isotropic Gaussian form

\[
M_{\alpha\beta\gamma\delta,\mathrm{eq}}^{(4)}
=
\rho T^2
\left(
\delta_{\alpha\beta}\delta_{\gamma\delta}
+
\delta_{\alpha\gamma}\delta_{\beta\delta}
+
\delta_{\alpha\delta}\delta_{\beta\gamma}
\right).
\]

Departures from these equilibrium central moments provide a useful physical interpretation of stress-like, heat-flux-like, and higher-order nonequilibrium structure. They are not, however, the quantities used directly by the implemented hierarchical sensor. The sensor and collision model described below use deviations of laboratory-frame raw Hermite coefficients from their local-equilibrium values.

\subsection{Hermite basis}

The fourth-order state used by the collision map is represented in a finite Hermite basis associated with the discrete quadrature \cite{grad1949kinetic,shan1998discretization}. Let

\[
\mathcal{H}^{(n)}
\left(
\boldsymbol{\xi}_i
\right)
\]

denote the rank-\(n\) probabilists' Hermite tensor associated with the standard-Gaussian weight and evaluated at the \(i\)-th discrete velocity. For example,

\[
\mathcal{H}^{(1)}_{\alpha}
=
\xi_{\alpha},
\qquad
\mathcal{H}^{(2)}_{\alpha\beta}
=
\xi_{\alpha}\xi_{\beta}
-
\delta_{\alpha\beta}.
\]

The truncated expansion through fourth order is

\[
f_i^{[4]}
=
w_i
\sum_{n=0}^{4}
\frac{1}{n!}
\boldsymbol{a}^{(n)}
:
\mathcal{H}^{(n)}
\left(
\boldsymbol{\xi}_i
\right),
\]

where \(\boldsymbol{a}^{(n)}\) is the rank-\(n\) coefficient tensor and the colon denotes complete contraction over all tensor indices.

More explicitly,

\[
f_i^{[4]}
=
w_i
\left[
a^{(0)}
+
a_{\alpha}^{(1)}
\mathcal{H}_{\alpha}^{(1)}
+
\frac{1}{2}
a_{\alpha\beta}^{(2)}
\mathcal{H}_{\alpha\beta}^{(2)}
+
\frac{1}{6}
a_{\alpha\beta\gamma}^{(3)}
\mathcal{H}_{\alpha\beta\gamma}^{(3)}
+
\frac{1}{24}
a_{\alpha\beta\gamma\delta}^{(4)}
\mathcal{H}_{\alpha\beta\gamma\delta}^{(4)}
\right].
\]

Repeated Cartesian indices are summed.

The numerical implementation constructs the basis tensors once for the D3Q125 velocity set and subsequently reuses them for projection, equilibrium reconstruction, nonequilibrium decomposition, and collision.

\subsection{Projection and reconstruction}

The Hermite coefficients are obtained by discrete projection,

\[
\boldsymbol{a}^{(n)}
=
\sum_i
f_i
\mathcal{H}^{(n)}
\left(
\boldsymbol{\xi}_i
\right).
\]

With the standard-Gaussian probabilists' convention, reconstruction uses the reciprocal factorial \(1/n!\), as displayed above. The five-node quadrature integrates products of retained Hermite tensors through fourth order exactly. Projection followed by reconstruction is therefore an identity on the retained completely symmetric \(H_0\)--\(H_4\) subspace, up to floating-point roundoff.

Given coefficient tensors through fourth order, the reconstructed distribution is denoted schematically by

\[
\mathcal{R}
\left[
\boldsymbol{a}^{(0)},
\boldsymbol{a}^{(1)},
\boldsymbol{a}^{(2)},
\boldsymbol{a}^{(3)},
\boldsymbol{a}^{(4)}
\right].
\]

Projection followed by reconstruction provides the finite-dimensional map used by the collision operator,

\[
f
\longrightarrow
\left\{
\boldsymbol{a}^{(n)}
\right\}_{n=0}^{4}
\longrightarrow
f^{*}.
\]

The zeroth- and first-order coefficients are not relaxed. Within the second-order coefficient, the trace associated with total energy is preserved, while only the deviatoric nonequilibrium component is relaxed. The third- and fourth-order nonequilibrium sectors are updated by their order-dependent relaxation factors, as described in Section 3.

\subsection{Local equilibrium}

The local equilibrium distribution is represented by the fourth-order Hermite projection of a Maxwellian with density \(\rho\), velocity \(\boldsymbol{u}\), and temperature \(T\),

\[
f_i^{\mathrm{eq}}
=
\mathcal{R}
\left[
\boldsymbol{a}_{\mathrm{eq}}^{(0)},
\boldsymbol{a}_{\mathrm{eq}}^{(1)},
\boldsymbol{a}_{\mathrm{eq}}^{(2)},
\boldsymbol{a}_{\mathrm{eq}}^{(3)},
\boldsymbol{a}_{\mathrm{eq}}^{(4)}
\right]_i.
\]

The equilibrium coefficient tensors are generated analytically from the local macroscopic fields and are evaluated in batch form for all spatial cells.

The equilibrium state satisfies the discrete conservation constraints

\[
\sum_i f_i^{\mathrm{eq}}
=
\rho,
\]

\[
\sum_i
\boldsymbol{\xi}_i
f_i^{\mathrm{eq}}
=
\rho\boldsymbol{u},
\]

and

\[
\sum_i
\frac{
\left|
\boldsymbol{\xi}_i
\right|^2
}{2}
f_i^{\mathrm{eq}}
=
E,
\]

where the local total discrete kinetic-energy density is

\[
E
=
\frac{1}{2}
\sum_i
\left|
\boldsymbol{\xi}_i
\right|^2
f_i
=
\frac{1}{2}
\rho
\left|
\boldsymbol{u}
\right|^2
+
\frac{3}{2}
\rho T.
\]

These identities hold up to floating-point roundoff for the stated fourth-order Hermite representation and D3Q125 quadrature.

\subsection{Nonequilibrium decomposition by order}

The projected coefficient tensors are decomposed as

\[
\boldsymbol{a}^{(n)}
=
\boldsymbol{a}_{\mathrm{eq}}^{(n)}
+
\boldsymbol{a}_{\mathrm{neq}}^{(n)}.
\]

For the retained nonconserved orders,

\[
\boldsymbol{a}_{\mathrm{neq}}^{(n)}
=
\boldsymbol{a}^{(n)}
-
\boldsymbol{a}_{\mathrm{eq}}^{(n)},
\qquad
n=2,3,4.
\]

The associated distribution-level contributions are

\[
f_i^{(n)}
=
\frac{w_i}{n!}
\boldsymbol{a}_{\mathrm{neq}}^{(n)}
:
\mathcal{H}^{(n)}
\left(
\boldsymbol{\xi}_i
\right).
\]

Consequently, the fourth-order Hermite reconstruction used by the collision map may be expressed as

\[
f_i^{[4]}
=
f_i^{\mathrm{eq}}
+
f_i^{(2)}
+
f_i^{(3)}
+
f_i^{(4)}.
\]

A population vector produced by the finite-volume transport stage need not lie exactly in this retained Hermite subspace. Projection followed by reconstruction therefore defines the \(H_0\)--\(H_4\) state used by the collision map and filters any component outside the retained representation. Because the projection retains the zeroth- and first-order coefficients and the complete second-order coefficient, the represented density, momentum, and total discrete kinetic energy are preserved up to quadrature and floating-point error.

This retained-subspace decomposition supplies both the order-specific TNE indicators and the independently relaxed collision channels.

\subsection{Scope of the fixed velocity set}

The present article keeps the D3Q125 velocity set fixed throughout the simulation. Adaptation is introduced only through the local order-specific relaxation factors; the velocity nodes, quadrature weights, and number of active populations do not change in space or time.

This distinction is important because it separates the contribution of hierarchical moment activation from the separate problem of adaptive velocity-space refinement. Dynamic switching among D3Q125, D3Q343, and D3Q729, velocity scaling, active-node selection, and conservative velocity-space remapping are deferred to future work.

\section{Hierarchical log-Gaussian collision model}

\subsection{Equilibrium and nonequilibrium decomposition}

Let \(f_i(\boldsymbol{x},t)\) denote the discrete population associated with the velocity \(\boldsymbol{\xi}_i\), \(i=1,\ldots,125\). At each spatial point, the collision map first projects the current population vector onto the retained fourth-order Hermite representation. The resulting reconstructed state is decomposed into an equilibrium contribution and retained nonequilibrium Hermite sectors,

\[
f_i^{[4]}
=
f_i^{\mathrm{eq}}
+
f_i^{(2)}
+
f_i^{(3)}
+
f_i^{(4)}.
\]

A transported population vector need not lie exactly in this subspace. The projection--reconstruction step therefore filters components outside the retained \(H_0\)--\(H_4\) representation before relaxation.

Here \(f_i^{\mathrm{eq}}\) is reconstructed from the local density, velocity, and temperature, while \(f_i^{(n)}\) denotes the contribution associated with the \(n\)-th-order nonequilibrium Hermite coefficients. The present formulation retains moments up to fourth order on the fixed D3Q125 velocity set.

For each order \(n\in\{2,3,4\}\), the nonequilibrium coefficient tensor is written schematically as

\[
\boldsymbol{a}^{(n)}_{\mathrm{neq}}
=
\boldsymbol{a}^{(n)}
-
\boldsymbol{a}^{(n)}_{\mathrm{eq}}.
\]

The collision operator acts independently on these retained sectors, while preserving the conserved density, momentum, and total-energy moments.

\subsection{Order-specific thermodynamic nonequilibrium measures}

The thermodynamic nonequilibrium sensor is evaluated from deviations of the projected raw Hermite coefficients from their local-equilibrium values. For orders \(n=2,3,4\),

\[
\Delta\boldsymbol{a}^{(n)}
=
\boldsymbol{a}^{(n)}
-
\boldsymbol{a}_{\mathrm{eq}}^{(n)}.
\]

The tensor magnitude is the Frobenius norm over all entries of the complete Cartesian tensor,

\[
\left\|
\Delta\boldsymbol{a}^{(n)}
\right\|_F
=
\left[
\sum_{\alpha_1,\ldots,\alpha_n}
\left(
\Delta a_{\alpha_1\cdots\alpha_n}^{(n)}
\right)^2
\right]^{1/2}.
\]

The implementation sums over the complete Cartesian representation. Consequently, symmetry-related permutations occur with their natural multiplicities. The order labels used by the sensor refer to raw Hermite sectors in the fixed laboratory frame; they are not claimed to define a strictly Galilean-invariant central-moment decomposition.

Let

\[
p_{\mathrm{loc}}=\rho T
\]

denote the local pressure. With the numerical floor

\[
\varepsilon=10^{-14},
\]

the pressure and temperature scales used in the implementation are

\[
\widehat p
=
\max(p_{\mathrm{loc}},\varepsilon),
\qquad
\widehat T
=
\max(T,\varepsilon).
\]

The three order-specific contributions are

\[
\mathcal E_2
=
c_2
\frac{
\left\|
\Delta\boldsymbol{a}^{(2)}
\right\|_F
}{
\widehat p
},
\]

\[
\mathcal E_3
=
c_3
\frac{
\left\|
\Delta\boldsymbol{a}^{(3)}
\right\|_F
}{
\widehat p\sqrt{\widehat T}
},
\]

and

\[
\mathcal E_4
=
c_4
\frac{
\left\|
\Delta\boldsymbol{a}^{(4)}
\right\|_F
}{
\widehat p\,\widehat T
}.
\]

Here \(c_2,c_3,c_4\) are configurable nonnegative coefficients. All numerical results reported in this article use

\[
c_2=c_3=c_4=1.
\]

The coefficients multiply the complete-Cartesian raw-Hermite deviation norms defined above.

The total TNE indicator is the arithmetic sum,

\[
\mathcal E_{\mathrm{tot}}
=
\mathcal E_2
+
\mathcal E_3
+
\mathcal E_4.
\]

It is not a \(p\)-norm. This distinction is important because the common sensor aggregates the three TNE channels additively before constructing the effective rarefaction indicator.

\subsection{Macroscopic-gradient rarefaction indicators}

For the one-dimensional periodic benchmarks, scalar gradients are computed with the second-order centered difference

\[
\left(
\partial_x q
\right)_j
=
\frac{
q_{j+1}-q_{j-1}
}{
2\Delta x
},
\]

with periodic indexing.

The density and temperature gradient magnitudes supplied to the sensor are

\[
G_\rho
=
\left|
\partial_x\rho
\right|,
\qquad
G_T
=
\left|
\partial_xT
\right|.
\]

For the three-component velocity field,

\[
\boldsymbol u
=
(u_x,u_y,u_z),
\]

the implementation applies the same centered difference independently to each component and then forms the Euclidean norm

\[
G_u
=
\left\|
\partial_x\boldsymbol u
\right\|_2
=
\left[
\left(\partial_xu_x\right)^2
+
\left(\partial_xu_y\right)^2
+
\left(\partial_xu_z\right)^2
\right]^{1/2}.
\]

The parameter \(\lambda\) is used here as a prescribed macroscopic-gradient sensor length scale. It multiplies the normalized gradients, but it does not independently prescribe a molecular collision time or the continuous-time collision frequency.

The macroscopic-gradient rarefaction indicators are

\[
K_\rho
=
\lambda
\frac{
G_\rho
}{
\max(\rho,\varepsilon)
},
\]

\[
K_T
=
\lambda
\frac{
G_T
}{
\max(T,\varepsilon)
},
\]

and

\[
K_u
=
\lambda
\frac{
G_u
}{
\max(c_{\mathrm{ref}},\varepsilon)
}.
\]

The velocity reference scale is

\[
c_{\mathrm{ref}}
=
\sqrt{
\max(T,\varepsilon)
}.
\]

Thus, the implementation uses the thermal scale \(\sqrt{T}\), not \(\sqrt{\gamma T}\).

The default numerical floor is

\[
\varepsilon=10^{-14}.
\]

When \(\lambda=0\), all three macroscopic-gradient indicators vanish exactly, while the TNE contribution remains active. This case is therefore called the \textbf{TNE-only sensor limit}, rather than a physical zero-mean-free-path limit.

\subsection{Order-resolved effective Knudsen indicators}

For nonnegative inputs \(q_1,\ldots,q_4\), define the unnormalized \(p\)-norm

\[
\mathcal N_p(q_1,q_2,q_3,q_4)
=
\left(
q_1^p+q_2^p+q_3^p+q_4^p
\right)^{1/p}.
\]

The default implementation uses

\[
p=8.
\]

For \(p=\infty\), the implementation instead returns the maximum input.

The order-resolved effective Knudsen indicators are

\[
K_n
=
\mathcal N_p
\left(
K_\rho,
K_T,
K_u,
\mathcal E_n
\right),
\qquad n=2,3,4.
\]

The common-sensor formulation uses

\[
K_{\mathrm{common}}
=
\mathcal N_p
\left(
K_\rho,
K_T,
K_u,
\mathcal E_{\mathrm{tot}}
\right)
\]

and supplies this same value to all retained moment orders.

Because the norm is not divided by the number of components, it is an unnormalized \(p\)-norm rather than a generalized mean.

\subsection{Order-dependent log-Gaussian relaxation spectrum}

The log-Gaussian schedule used here follows the scale-space weighting proposed in \cite{wu2026scaleSpace}. The present construction applies that probability partition separately to retained Hermite orders rather than to a blended interface flux.

For each moment order, the continuum fraction is

\[
w_{c,n}(K)
=
\frac{1}{2}
\operatorname{erfc}
\left[
\frac{
\ln\!\left(
\widehat K/K_{0,n}
\right)
}{
\sqrt{2}\,\sigma_n
}
\right],
\]

where

\[
\widehat K
=
\max(K,K_{\mathrm{floor}})
\]

and the default floor is

\[
K_{\mathrm{floor}}=10^{-14}.
\]

The kinetic fraction is

\[
w_{f,n}(K)=1-w_{c,n}(K).
\]

The calibrated reference-step relaxation factor is the convex blend

\[
s_{n,\mathrm{ref}}(K)
=
w_{c,n}(K)s_{n,\mathrm{cont}}
+
w_{f,n}(K)s_{n,\mathrm{kin}}.
\]

The stored endpoint values are therefore discrete factors calibrated at the reference timestep

\[
\Delta t_{\mathrm{ref}}
=
1.0938161705673147\times 10^{-3}.
\]

For an arbitrary positive timestep, the timestep-consistent discrete factor is defined directly through the survival mapping

\[
1-s_n(K_n,\Delta t)
=
\left[
1-s_{n,\mathrm{ref}}(K_n)
\right]^{
\Delta t/\Delta t_{\mathrm{ref}}
},
\]

or equivalently,

\[
s_n(K_n,\Delta t)
=
1-
\left[
1-s_{n,\mathrm{ref}}(K_n)
\right]^{
\Delta t/\Delta t_{\mathrm{ref}}
}.
\]

For \(0\leq s_{n,\mathrm{ref}}(K_n)<1\), the same mapping may be written using the finite continuous-time frequency

\[
\nu_n(K_n)
=
-\frac{
\ln\!\left[
1-s_{n,\mathrm{ref}}(K_n)
\right]
}{
\Delta t_{\mathrm{ref}}
},
\]

so that

\[
s_n(K_n,\Delta t)
=
1-\exp\!\left[
-\nu_n(K_n)\Delta t
\right].
\]

When \(s_{n,\mathrm{ref}}=1\), the survival mapping is used directly. For every positive timestep, its survival factor is zero and hence \(s_n=1\). This is the complete-relaxation limit, corresponding formally to \(\nu_n\rightarrow+\infty\).

Consequently, two half-step collision maps have exactly the same survival factor as one full-step collision map when \(K_n\) is held fixed.

The actual-timestep factor is applied in the local collision map,

\[
\boldsymbol{a}_{\mathrm{neq}}^{(n),*}
=
\left[
1-s_n(K_n,\Delta t)
\right]
\boldsymbol{a}_{\mathrm{neq}}^{(n)}.
\]

The endpoint values in the spectrum are therefore reference-step calibration factors rather than continuous-time frequencies. The exponential survival mapping makes the accumulated collision action consistent under timestep subdivision when the sensor value \(K_n\) is held fixed.

The default parameters are

\[
\begin{array}{c|cccc}
n
&
K_{0,n}
&
\sigma_n
&
s_{n,\mathrm{cont}}
&
s_{n,\mathrm{kin}}
\\
\hline
2 & 0.050 & 2.0 & 1.0 & 0.20 \\
3 & 0.030 & 2.5 & 1.0 & 0.10 \\
4 & 0.015 & 3.0 & 1.0 & 0.05
\end{array}.
\]

No additional clipping of the factor is performed. For the default parameters, the erfc weight keeps each factor between its continuum and kinetic limits.

\subsection{Raw-Hermite collision update and positivity control}

The implemented collision acts on raw Hermite coefficient tensors, not directly on central moments.

The zeroth- and first-order coefficients are copied unchanged,

\[
a^{(0),*}=a^{(0)},
\qquad
\boldsymbol{a}^{(1),*}
=
\boldsymbol{a}^{(1)}.
\]

For the second-order coefficient, define the pre-collision trace

\[
\tau
=
\operatorname{tr}
\boldsymbol{a}^{(2)}
\]

and the associated isotropic tensor

\[
\boldsymbol{a}^{(2)}_{\mathrm{iso}}
=
\frac{\tau}{3}\boldsymbol{I}.
\]

The pre-collision deviatoric part is

\[
\boldsymbol{a}^{(2)}_{\mathrm{dev}}
=
\boldsymbol{a}^{(2)}
-
\boldsymbol{a}^{(2)}_{\mathrm{iso}}.
\]

Similarly, the equilibrium deviatoric part is

\[
\boldsymbol{a}^{(2)}_{\mathrm{eq,dev}}
=
\boldsymbol{a}^{(2)}_{\mathrm{eq}}
-
\frac{
\operatorname{tr}
\boldsymbol{a}^{(2)}_{\mathrm{eq}}
}{3}
\boldsymbol{I}.
\]

The post-collision second-order tensor is

\[
\boldsymbol{a}^{(2),*}
=
\boldsymbol{a}^{(2)}_{\mathrm{iso}}
+
\boldsymbol{a}^{(2)}_{\mathrm{eq,dev}}
+
(1-s_2)
\left(
\boldsymbol{a}^{(2)}_{\mathrm{dev}}
-
\boldsymbol{a}^{(2)}_{\mathrm{eq,dev}}
\right).
\]

Thus, the pre-collision second-order trace is preserved and only the deviatoric component is relaxed.

For the standard-Gaussian Hermite basis, the local discrete kinetic energy can be written as

\[
E
=
\frac{1}{2}
\left[
3a^{(0)}
+
\operatorname{tr}
\boldsymbol{a}^{(2)}
\right].
\]

Because \(a^{(0)}\) and the pre-collision second-order trace are both preserved, the local discrete kinetic energy is unchanged by collision up to floating-point roundoff.

The third- and fourth-order tensors use the standard linear relaxation

\[
\boldsymbol{a}^{(3),*}
=
\boldsymbol{a}^{(3)}_{\mathrm{eq}}
+
(1-s_3)
\left(
\boldsymbol{a}^{(3)}
-
\boldsymbol{a}^{(3)}_{\mathrm{eq}}
\right),
\]

\[
\boldsymbol{a}^{(4),*}
=
\boldsymbol{a}^{(4)}_{\mathrm{eq}}
+
(1-s_4)
\left(
\boldsymbol{a}^{(4)}
-
\boldsymbol{a}^{(4)}_{\mathrm{eq}}
\right).
\]

The post-collision population is reconstructed from the updated raw Hermite coefficients.

When positivity control is requested, let \(f_i^{\mathrm{cand}}\) denote the unconstrained candidate and \(f_i^{\mathrm{eq}}\) the reconstructed truncated-Hermite equilibrium distribution corresponding to the same local macroscopic state. The limited distribution is

\[
f_i^{\mathrm{lim}}
=
f_i^{\mathrm{eq}}
+
\theta
\left(
f_i^{\mathrm{cand}}
-
f_i^{\mathrm{eq}}
\right),
\qquad
0\le\theta\le1.
\]

For every violating population satisfying

\[
f_i^{\mathrm{cand}}
<
f_{\mathrm{floor}}-\varepsilon_{\mathrm{tol}},
\]

the admissible population factor is

\[
\theta_i
=
\frac{
f_i^{\mathrm{eq}}
-
f_{\mathrm{floor}}
}{
f_i^{\mathrm{eq}}
-
f_i^{\mathrm{cand}}
}.
\]

The cellwise blending factor is

\[
\theta
=
\operatorname{clip}
\left(
\min_i\theta_i,
0,
1
\right),
\]

with nonviolating populations contributing the value \(1\).

The equilibrium reference must satisfy

\[
f_i^{\mathrm{eq}}
\ge
f_{\mathrm{floor}}
-
\varepsilon_{\mathrm{tol}}
\]

for every velocity; otherwise the limiter raises an error. The default tolerance is

\[
\varepsilon_{\mathrm{tol}}=10^{-14}.
\]

If every candidate population is already at or above the requested floor, then \(\theta=1\) and the candidate is returned unchanged. Populations lying below the floor but within the roundoff tolerance do not trigger convex blending; they are instead snapped componentwise to the floor.

Because both endpoints of the line segment are reconstructed from the same local conserved state, the convex-blending stage preserves their shared conserved moments up to reconstruction and floating-point error. The final componentwise snap of populations lying within \(\varepsilon_{\mathrm{tol}}\) below the floor can introduce an additional roundoff-scale perturbation in mass, momentum, and energy. A cell is reported as limited whenever either convex blending or this final snap modifies its population vector.

\subsection{Conservation and entropy diagnostics}

The collision and transport stages are monitored separately. The reported diagnostics include changes in total mass,

\[
\Delta M
=
\sum_{\boldsymbol{x},i}
f_i^{\mathrm{after}}
-
\sum_{\boldsymbol{x},i}
f_i^{\mathrm{before}},
\]

total momentum,

\[
\Delta\boldsymbol{P}
=
\sum_{\boldsymbol{x},i}
\boldsymbol{\xi}_i
\left(
    f_i^{\mathrm{after}}
    -
    f_i^{\mathrm{before}}
\right),
\]

and kinetic energy,

\[
\Delta E
=
\frac{1}{2}
\sum_{\boldsymbol{x},i}
\left|\boldsymbol{\xi}_i\right|^2
\left(
    f_i^{\mathrm{after}}
    -
    f_i^{\mathrm{before}}
\right).
\]

For positive distributions, the discrete Boltzmann functional

\[
H[f]
=
\sum_{\boldsymbol{x},i}
f_i
\ln
\left(
    \frac{f_i}{w_i}
\right)
\]

is also recorded. The present work uses these quantities as numerical diagnostics rather than claiming a general discrete \(H\)-theorem for the truncated hierarchical collision operator. Entropic lattice Boltzmann formulations instead construct the collision step from a discrete entropy and enforce an associated entropy condition \cite{succi2002htheorem,ansumali2002entropic}; that construction is not used here.

\subsection{Limiting behavior}

The order-resolved model has two important limiting regimes.

When the macroscopic-gradient contribution dominates,

\[
K_\rho,\ K_T,\ K_u
\gg
\mathcal{E}_2,\mathcal{E}_3,\mathcal{E}_4,
\]

the effective indicators satisfy approximately

\[
K_2\approx K_3\approx K_4,
\]

and the hierarchical model approaches the common-sensor response.

When the shared macroscopic background is weak and the nonequilibrium channels differ appreciably,

\[
\mathcal{E}_2
\ne
\mathcal{E}_3
\ne
\mathcal{E}_4,
\]

the three relaxation factors separate. This is the regime in which the order-resolved mechanism has its largest dynamical effect.

\section{Spatial discretization and time integration}

\subsection{Finite-volume transport equation}

Finite-volume discretizations provide an off-lattice route for advancing discrete kinetic populations on spatial meshes that need not coincide with the discrete-velocity lattice \cite{peng1999finitevolume,patil2009tvd}. For the one-dimensional benchmarks considered in this work, the discrete kinetic equation is written as

\[
\frac{\partial f_i}{\partial t}
+
\xi_{i,x}
\frac{\partial f_i}{\partial x}
=
\Omega_i.
\]

The spatial domain is divided into uniform finite-volume cells of width

\[
\Delta x
=
\frac{L}{N_x},
\]

where \(L\) is the domain length and \(N_x\) is the number of cells. The cell-averaged distribution is denoted by

\[
f_{i,j}(t),
\]

with \(j\) the spatial-cell index.

The transport and collision stages are applied sequentially within each timestep. The collision operator acts locally in each cell, whereas the transport operator communicates populations between neighboring cells according to the sign of the streamwise discrete velocity \(\xi_{i,x}\).

\subsection{First-order upwind transport}

The first-order upwind discretization is used as the baseline transport scheme. The semi-discrete finite-volume update is

\[
\frac{\mathrm{d}f_{i,j}}{\mathrm{d}t}
=
-
\frac{
F_{i,j+1/2}
-
F_{i,j-1/2}
}{
\Delta x
},
\]

where the numerical flux is

\[
F_{i,j+1/2}
=
\xi_{i,x}
f_{i,j+1/2}^{\mathrm{up}}.
\]

The upwind interface state is

\[
f_{i,j+1/2}^{\mathrm{up}}
=
\begin{cases}
f_{i,j},
&
\xi_{i,x}\ge 0,
\\[4pt]
f_{i,j+1},
&
\xi_{i,x}<0.
\end{cases}
\]

The corresponding forward-Euler transport update is

\[
f_{i,j}^{n+1}
=
f_{i,j}^{*}
-
\frac{
\Delta t
}{
\Delta x
}
\left(
F_{i,j+1/2}^{*}
-
F_{i,j-1/2}^{*}
\right),
\]

where \(f_{i,j}^{*}\) denotes the post-collision population.

The first-order scheme is robust and positivity friendly under the standard CFL restriction, but it introduces numerical diffusion. It is therefore used primarily for mechanism isolation and regression tests.

\subsection{MUSCL reconstruction}

A second-order MUSCL finite-volume option is also available. Within each cell, a limited slope is reconstructed from neighboring cell averages. For a scalar population component, the minmod limiter is written as

\[
\operatorname{minmod}(a,b)
=
\begin{cases}
\operatorname{sign}(a)
\min(|a|,|b|),
&
ab>0,
\\[4pt]
0,
&
ab\le 0.
\end{cases}
\]

The limited slope in cell \(j\) is

\[
\sigma_{i,j}
=
\operatorname{minmod}
\left(
f_{i,j}-f_{i,j-1},
f_{i,j+1}-f_{i,j}
\right).
\]

The left and right reconstructed states are

\[
f_{i,j}^{L}
=
f_{i,j}
-
\frac{1}{2}
\sigma_{i,j},
\]

and

\[
f_{i,j}^{R}
=
f_{i,j}
+
\frac{1}{2}
\sigma_{i,j}.
\]

At the interface \(j+1/2\), the upwind flux uses

\[
f_{i,j+1/2}^{\mathrm{MUSCL}}
=
\begin{cases}
f_{i,j}^{R},
&
\xi_{i,x}\ge 0,
\\[4pt]
f_{i,j+1}^{L},
&
\xi_{i,x}<0.
\end{cases}
\]

The limiter suppresses spurious oscillations while retaining second-order accuracy in smooth regions. Related finite-volume lattice kinetic formulations use TVD reconstruction for the physical-space advection of discrete populations \cite{patil2009tvd}.

\subsection{Time integration}

The timestep is restricted by the maximum streamwise discrete speed,

\[
\Delta t
\le
\mathrm{CFL}
\frac{
\Delta x
}{
\max_i |\xi_{i,x}|
}.
\]

The benchmarks use a prescribed CFL number and compute the timestep directly from the velocity set and spatial resolution.

Each timestep is implemented as a sequential collision--transport update. The adaptive sensor and the order-dependent relaxation factors are evaluated from the distribution at the beginning of the step. One local collision map is then applied, followed by either the first-order upwind transport routine or the MUSCL transport routine.

The adaptive sensor and collision operator are not reevaluated within transport substages. The algorithm should therefore be interpreted as a discrete collision-then-transport splitting, not as an SSPRK2 method-of-lines integration of a combined collision--transport operator.

The same timestep is used for all discrete populations.

\subsection{Automatic gradient evaluation for the adaptive sensor}

For periodic one-dimensional states, the adaptive rarefaction sensor uses the centered finite difference

\[
\left(
\partial_x q
\right)_j
=
\frac{
q_{j+1}-q_{j-1}
}{
2\Delta x
}.
\]

The density and temperature channels use the absolute values of these derivatives. For velocity, the derivative is evaluated componentwise for all three velocity components and combined through the Euclidean norm,

\[
\left\|
\partial_x\boldsymbol u
\right\|_2.
\]

The velocity-gradient channel is normalized by

\[
\sqrt{\max(T,10^{-14})}.
\]

These sensor gradients are independent of the upwind or MUSCL reconstruction used for population transport.

\subsection{Collision-transport sequence}

At each timestep, the following operations are performed:

\begin{enumerate}
\item evaluate density, velocity, and temperature, and project the retained raw Hermite coefficients from the current distribution;
\item compute the macroscopic-gradient rarefaction indicators;
\item compute the second-, third-, and fourth-order TNE indicators;
\item construct either one common effective Knudsen indicator or the three order-resolved indicators \(K_2,K_3,K_4\);
\item evaluate the order-dependent relaxation factors;
\item apply the local moment-space collision operator;
\item apply the optional positivity limiter;
\item transport the post-collision distribution;
\item evaluate invariant, positivity, and entropy diagnostics.
\end{enumerate}

This sequence keeps the sensor evaluation and collision response local, while the finite-volume stage accounts for spatial propagation.

\subsection{Periodic boundary conditions}

The smooth temperature-wave and compression-wave benchmarks use periodic boundaries. For any discrete population,

\[
f_{i,0}
=
f_{i,N_x},
\]

and

\[
f_{i,N_x+1}
=
f_{i,1},
\]

in ghost-cell notation.

Periodic boundaries eliminate net boundary fluxes, so changes in total mass, momentum, and energy can be attributed to the numerical collision-transport procedure rather than external forcing.

\subsection{Positivity monitoring}

The minimum population is recorded at each diagnostic step,

\[
f_{\min}
=
\min_{\boldsymbol{x},i}
f_i.
\]

A positive value confirms that the discrete entropy functional is defined and that the state remains admissible as a nonnegative discrete-population representation.

When the optional positivity limiter is active, it is applied after collision and before transport. The numerical benchmarks reported here also record the unconstrained minimum population to determine whether limiting is actually required.

\subsection{Conservation diagnostics}

The total mass is

\[
M
=
\sum_j\sum_i
f_{i,j}.
\]

The total momentum is

\[
\boldsymbol{P}
=
\sum_j\sum_i
\boldsymbol{\xi}_i f_{i,j},
\]

and the total kinetic energy is

\[
E
=
\frac{1}{2}
\sum_j\sum_i
|\boldsymbol{\xi}_i|^2
f_{i,j}.
\]

Because every reported spatial mesh is uniform, these quantities are stored and reported as discrete cell sums with the common cell-width factor \(\Delta x\) omitted. Restoring the corresponding spatial-integral convention would multiply \(M\), \(\boldsymbol{P}\), and \(E\), together with their absolute changes, by \(\Delta x\). This omission does not affect same-grid relative conservation errors or the periodic flux-telescoping identities. Absolute grid sums from different spatial resolutions should not be compared without restoring this factor.

The solver records collision-stage changes,

\[
\Delta M_{\mathrm{coll}},
\qquad
\Delta\boldsymbol{P}_{\mathrm{coll}},
\qquad
\Delta E_{\mathrm{coll}},
\]

and net boundary-transport changes,

\[
\Delta M_{\mathrm{boundary}},
\qquad
\Delta\boldsymbol{P}_{\mathrm{boundary}},
\qquad
\Delta E_{\mathrm{boundary}}.
\]

For periodic boundaries, the conservative flux differences telescope globally for each discrete population. The resulting transport-stage global changes therefore vanish up to floating-point accumulation error. Spatial discretization affects the accuracy of the transported solution but not this global conservation identity.

\subsection{Entropy diagnostics}

For strictly positive populations, the discrete entropy-like functional is evaluated as

\[
H[f]
=
\sum_j\sum_i
f_{i,j}
\ln
\left(
\frac{f_{i,j}}{w_i}
\right).
\]

As for the invariant diagnostics, the common cell-width factor \(\Delta x\) is omitted from this uniform-grid sum. Including it would multiply \(H\), its stagewise changes, and the entropy-accounting residual by the same constant factor. It would not alter the bookkeeping-closure identity or any sign statement made for a fixed grid.

The code records the collision-stage change,

\[
\Delta H_{\mathrm{coll}},
\]

and transport-stage change,

\[
\Delta H_{\mathrm{trans}}.
\]

The stagewise entropy-accounting residual is defined by comparing the total change with the sum of the collision- and transport-stage changes. It is a bookkeeping-closure diagnostic for the split update, not a physical entropy-production residual. These quantities are used to detect unexpected numerical behavior, but no general monotonicity theorem is assumed for every transport and limiting configuration.

\subsection{Fixed-final-time grid-and-timestep sensitivity}

The tested cell-count and timestep-count pairs are

\[
(N_x,N_t)=(32,1),(64,2),(128,4),(256,8),
\]

with

\[
\Delta t \propto \Delta x \propto N_x^{-1}.
\]

All four calculations therefore reach the same nominal final time,

\[
t_f=4.375265\times10^{-3}.
\]

The calibrated reference-step factors are mapped to each actual timestep through exponential survival scaling,

\[
1-s_n(K_n,\Delta t)
=
\left[
1-s_{n,\mathrm{ref}}(K_n)
\right]^{
\Delta t/\Delta t_{\mathrm{ref}}
}.
\]

For a fixed sensor value, timestep subdivision therefore preserves the accumulated collision survival factor. The calculations still vary the spatial resolution, timestep, sequential collision--transport splitting error, and dynamically evolving sensor fields together. The experiment is therefore interpreted as a fixed-final-time grid-and-timestep sensitivity study rather than a pure spatial-convergence test.

\section{Numerical experiments}

\subsection{Pure-order hierarchical activation}

The first set of tests isolates the response of the collision operator to nonequilibrium perturbations confined to a single retained moment order. Separate initial distributions are constructed with dominant second-, third-, or fourth-order Hermite content, while the remaining nonequilibrium sectors are kept at roundoff level.

For a pure order-\(n\) perturbation, the corresponding order-specific sensor \(K_n\) becomes the dominant effective Knudsen indicator. The other two indicators remain close to the common macroscopic background. Consequently, the matching relaxation factor responds most strongly to the imposed perturbation, while nonmatching channels remain essentially unchanged.

The measured matching responses are consistent with the prescribed order-dependent relaxation spectrum. Representative effective Knudsen indicators and matching relaxation factors are

\[
K_2 \approx 7.35\times10^{-2},
\qquad
s_2 \approx 5.39\times10^{-1},
\]

\[
K_3 \approx 3.87\times10^{-2},
\qquad
s_3 \approx 5.13\times10^{-1},
\]

and

\[
K_4 \approx 3.35\times10^{-2},
\qquad
s_4 \approx 4.25\times10^{-1}.
\]

For each pure-order test, the decay of the matching nonequilibrium coefficient agrees with the expected factor

\[
1-s_n,
\]

while the nonmatching signals remain at floating-point level. This confirms that the implementation preserves the intended diagonal moment-space structure and does not introduce artificial cross-order activation.

\begin{figure}[tbp]
\centering
\includegraphics[width=\textwidth,height=0.68\textheight,keepaspectratio]{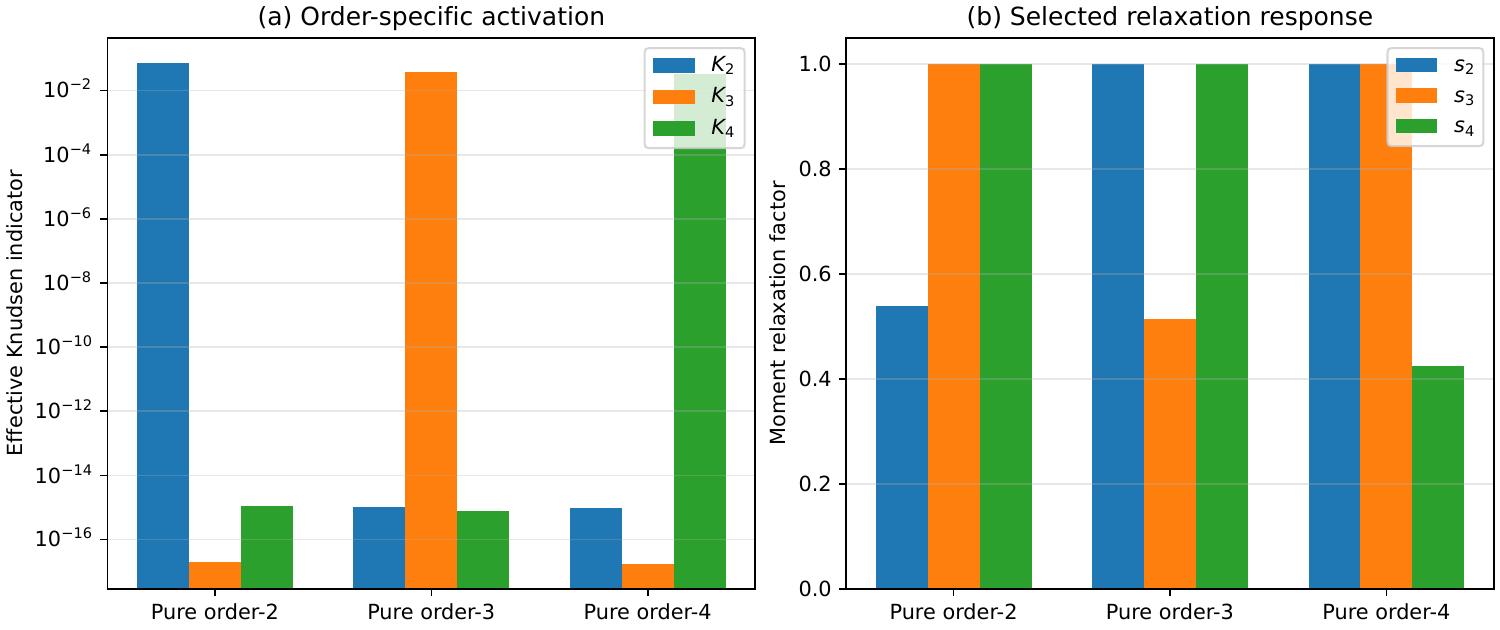}
\caption{Pure-order activation of the hierarchical relaxation spectrum. Separate second-, third-, and fourth-order perturbations activate their matching effective Knudsen indicators and relaxation factors. Nonmatching moment sectors remain negligible, demonstrating order selectivity of the collision response.}
\label{fig:pure-order-activation}
\end{figure}

\subsection{Homogeneous mixed-order relaxation}

The second test considers a spatially homogeneous distribution containing simultaneous second-, third-, and fourth-order nonequilibrium perturbations. Because transport is absent, the experiment isolates the cumulative effect of repeated collision updates.

Two formulations are compared. In the common-sensor model, the three retained moment sectors share a single effective Knudsen indicator constructed from the total nonequilibrium intensity. In the order-resolved model, the second-, third-, and fourth-order sectors use their own indicators \(K_2\), \(K_3\), and \(K_4\).

The common formulation couples all moment sectors through the largest contribution to the total sensor. A strongly activated channel can therefore suppress the relaxation of otherwise weakly nonequilibrium orders. The order-resolved model removes this artificial cross-order constraint and permits each sector to relax according to its local nonequilibrium magnitude.

Across the homogeneous relaxation history, the order-resolved model produces a smaller residual total TNE than the common model while maintaining positive populations. The entropy-like functional and minimum population remain well behaved throughout the collision sequence. The result demonstrates that hierarchical activation changes the relaxation pathway even in the absence of spatial gradients.

\subsection{Robustness to perturbation amplitude and composition}

\subsubsection{Amplitude scan}

To test whether the hierarchical effect depends on one specially chosen perturbation magnitude, the mixed-order initial state is scaled over a broad amplitude range. The relative second-, third-, and fourth-order composition is held fixed while the total nonequilibrium intensity is varied.

For all tested amplitudes, the residual total TNE produced by the order-resolved formulation remains below that of the common-sensor formulation. The ratio

\[
R_{\mathrm{TNE}}
=
\frac{
\mathcal{E}_{\mathrm{resolved}}
}{
\mathcal{E}_{\mathrm{common}}
}
\]

ranges approximately from

\[
0.556
\]

at the weakest perturbation to

\[
0.786
\]

at the strongest perturbation. Thus, the hierarchical advantage is largest in the weakly nonequilibrium regime but persists across the entire tested amplitude interval.

All order-specific relaxation factors remain positive, and the population minimum remains above zero. The amplitude scan therefore supports both the robustness and numerical admissibility of the hierarchical collision response.

\subsubsection{Composition scan}

The dependence on the distribution of nonequilibrium content among moment orders is examined using second-dominant, third-dominant, and fourth-dominant initial states. The prescribed relative amplitudes are

\[
(1.5,0.5,0.5),
\]

\[
(0.5,1.5,0.5),
\]

and

\[
(0.25,0.5,2.0)
\]

for the second-, third-, and fourth-order sectors, respectively.

In every case, the requested sector is the dominant initial TNE contributor. The order-resolved-to-common residual-TNE ratios are

\[
0.841
\]

for the second-dominant case,

\[
0.731
\]

for the third-dominant case, and

\[
0.750
\]

for the fourth-dominant case.

The hierarchical formulation therefore reduces the residual nonequilibrium response for all three order compositions. The magnitude of the improvement varies with the distribution of TNE among the retained sectors, but the sign of the effect is unchanged.

All tested distributions remain positive. These results show that the hierarchical mechanism is not tied to one dominant moment order and does not rely on a narrowly tuned initial composition.

\begin{figure}[tbp]
\centering
\includegraphics[width=\textwidth,height=0.68\textheight,keepaspectratio]{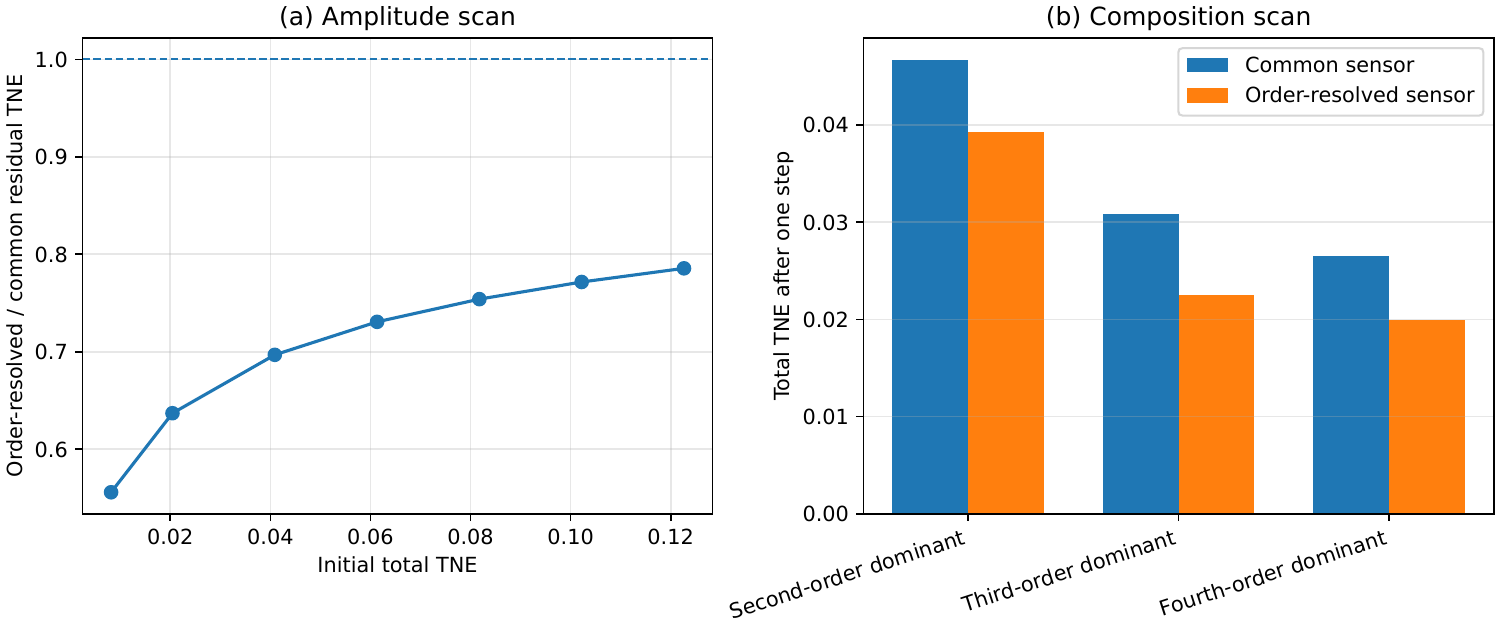}
\caption{Robustness of hierarchical activation to perturbation strength and order composition. (a) Ratio of order-resolved to common residual total TNE across the amplitude scan. The ratio remains below unity for all tested perturbation magnitudes. (b) Residual total TNE for second-, third-, and fourth-dominant initial states. The order-resolved formulation yields a lower residual in every composition.}
\label{fig:hierarchical-activation-summary}
\end{figure}

\subsection{Smooth temperature-wave benchmark}

\subsubsection{Numerical setup}

The first spatial benchmark considers a smooth periodic temperature wave on the fixed D3Q125 velocity set. The density is uniform, the macroscopic velocity is initially zero, and the temperature is prescribed as

\[
T(x,0)
=
1
+
0.15
\cos
\left(
\frac{2\pi x}{L}
\right).
\]

The initial distribution is reconstructed from the corresponding local equilibrium state. The domain is discretized with 128 cells, periodic boundary conditions are applied, and the macroscopic-gradient sensor scale is set to

\[
\lambda=0.01.
\]

The timestep is determined from the discrete-velocity CFL condition. Profiles are recorded at steps 0, 1, 2, and 4.

Because the initial density and velocity fields are spatially uniform, the temperature-gradient contribution is the dominant macroscopic rarefaction channel. At the initial state,

\[
K_\rho=0,
\qquad
K_u=0,
\]

up to numerical precision, while \(K_T\) has the imposed sinusoidal structure.

\subsubsection{Generation of order-specific nonequilibrium}

The initial local-equilibrium reconstruction produces second-, third-, and fourth-order TNE contributions at roundoff level. After transport begins, finite nonequilibrium is generated in all three retained moment sectors.

At step 4, the maximum recorded values are approximately

\[
\max \mathcal{E}_2
=
6.55\times10^{-5},
\]

\[
\max \mathcal{E}_3
=
5.60\times10^{-3},
\]

and

\[
\max \mathcal{E}_4
=
3.14\times10^{-4}.
\]

The third-order nonequilibrium contribution is therefore dominant by more than one order of magnitude. This response is consistent with the transport of a spatially varying thermal field, which generates a strong heat-flux-like third-order signal before comparable second- or fourth-order deviations develop.

The order-specific effective Knudsen indicators remain close to the shared temperature-gradient background, but the third-order channel acquires the largest TNE-dependent correction. At step 4,

\[
\max K_2
\approx
9.519\times10^{-3},
\]

\[
\max K_3
\approx
9.536\times10^{-3},
\]

and

\[
\max K_4
\approx
9.519\times10^{-3}.
\]

The small absolute separation reflects the fact that the imposed macroscopic temperature gradient is still the dominant component of the sensor at \(\lambda=0.01\).

\subsubsection{Hierarchical relaxation response}

Despite the small separation among the effective Knudsen indicators, the order-dependent relaxation spectra produce clearly distinct local factors. At step 4, the relaxation-factor ranges are approximately

\[
0.837
\le
s_2
\le
0.997,
\]

\[
0.709
\le
s_3
\le
0.976,
\]

and

\[
0.582
\le
s_4
\le
0.920.
\]

The hierarchy follows the prescribed order-dependent activation curves: the fourth-order sector admits the slowest relaxation, the third-order sector is intermediate, and the second-order sector remains closest to its continuum-limit factor.

The minimum discrete population remains positive,

\[
\min_i f_i
=
5.31\times10^{-7},
\]

so the benchmark does not require positivity limiting.

\begin{figure}[tbp]
\centering
\includegraphics[width=\textwidth,height=0.68\textheight,keepaspectratio]{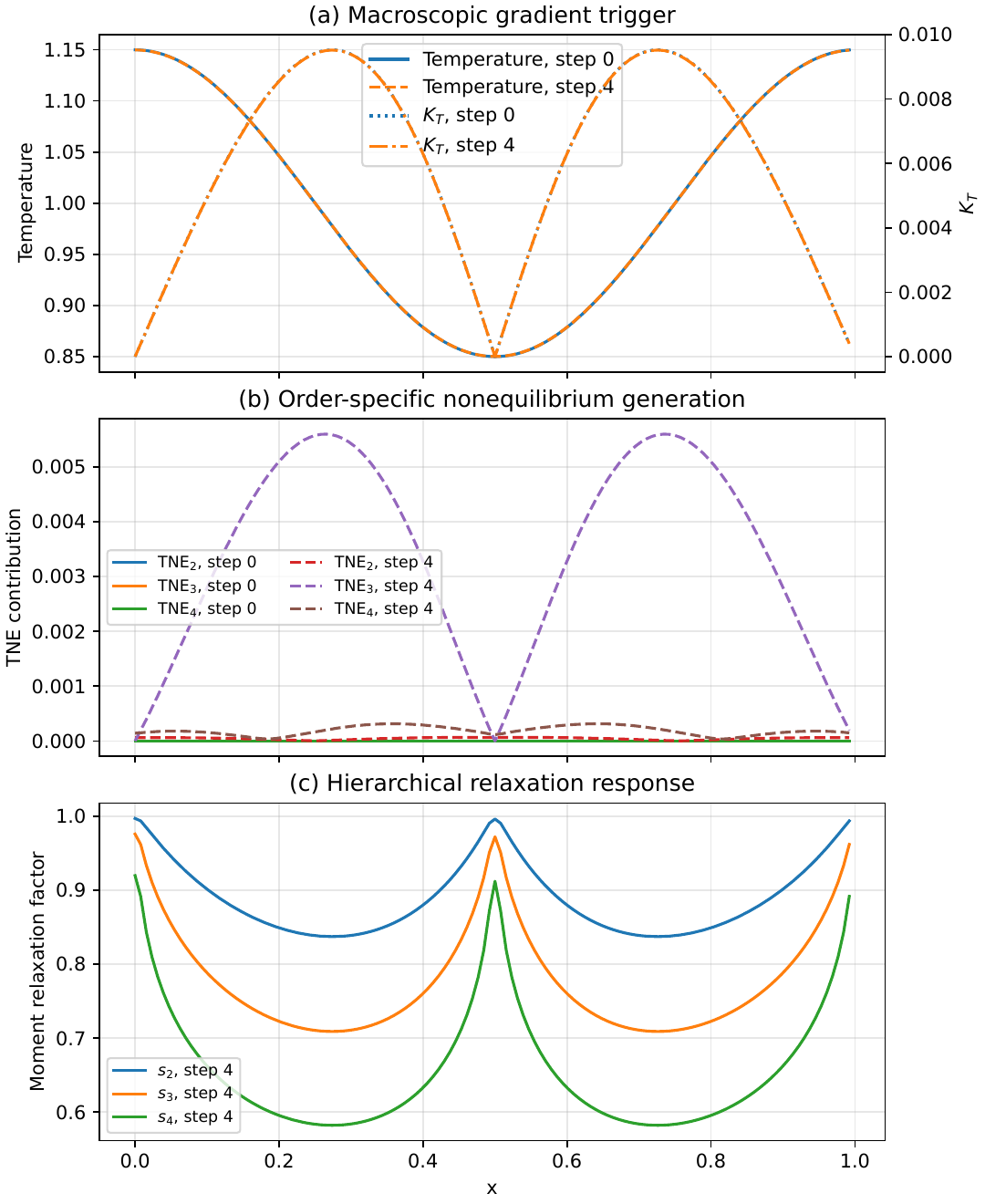}
\caption{Order-resolved activation in a smooth periodic temperature wave. (a) Temperature field and temperature-gradient rarefaction indicator at the initial state and after four timesteps. (b) Second-, third-, and fourth-order thermodynamic nonequilibrium contributions. The third-order signal dominates after transport generates a finite heat-flux-like response. (c) Corresponding order-dependent relaxation factors at step 4. The three moment sectors exhibit clearly separated collision responses despite sharing nearly the same macroscopic-gradient background.}
\label{fig:smooth-temperature-wave}
\end{figure}

\subsubsection{Interpretation}

The smooth temperature wave provides a controlled example in which one macroscopic-gradient channel dominates while order-specific TNE is generated dynamically. The initial state contains no meaningful nonequilibrium content, so the subsequent separation of \(\mathcal{E}_2\), \(\mathcal{E}_3\), and \(\mathcal{E}_4\) is caused by the transport-collision evolution rather than by an imposed kinetic perturbation.

The experiment demonstrates two complementary features of the model. First, the TNE sensor correctly identifies the third-order sector as the dominant nonequilibrium channel. Second, the distinct relaxation spectra preserve a clear hierarchy among \(s_2\), \(s_3\), and \(s_4\), even when the effective Knudsen indicators differ only weakly.

This benchmark is intentionally smooth and near equilibrium. It is therefore used to verify selective activation and positivity rather than to demonstrate a large difference between the common and order-resolved formulations. A stronger multi-gradient test is provided by the smooth compression wave in Section 5.5.

\subsection{Smooth compression-wave benchmark}

\subsubsection{Numerical setup}

The second spatial benchmark considers a smooth periodic compression wave in which density, streamwise velocity, and temperature vary simultaneously. The initial fields are

\[
\rho(x,0)
=
1
+
0.08
\cos
\left(
\frac{2\pi x}{L}
\right),
\]

\[
u_x(x,0)
=
0.06
\sin
\left(
\frac{2\pi x}{L}
\right),
\]

and

\[
T(x,0)
=
1
+
0.08
\cos
\left(
\frac{2\pi x}{L}
+
\frac{\pi}{4}
\right).
\]

The remaining velocity components are zero. The phase shift between the three fields ensures that the density-, temperature-, and velocity-gradient rarefaction indicators are simultaneously active but do not reach their extrema at the same spatial locations.

The initial distribution is reconstructed from the corresponding local equilibrium state on the fixed D3Q125 velocity set. Unless stated otherwise, the domain is discretized with 128 cells, periodic boundary conditions are used, the macroscopic-gradient sensor scale is

\[
\lambda=0.01,
\]

and profiles are recorded through four timesteps.

\subsubsection{Joint macroscopic-gradient activation}

At the initial state, the maximum macroscopic rarefaction contributions are approximately

\[
\max K_\rho
=
5.04\times10^{-3},
\]

\[
\max K_T
=
5.04\times10^{-3},
\]

and

\[
\max K_u
=
3.88\times10^{-3}.
\]

The second-, third-, and fourth-order TNE contributions are initially at roundoff level. Transport subsequently generates finite nonequilibrium signals in all three retained moment sectors.

After four timesteps,

\[
\max \mathcal{E}_2
=
7.98\times10^{-4},
\]

\[
\max \mathcal{E}_3
=
2.73\times10^{-3},
\]

and

\[
\max \mathcal{E}_4
=
8.28\times10^{-4}.
\]

The third-order sector is dominant, while the second- and fourth-order signals remain finite and spatially structured. The corresponding effective Knudsen indicators are

\[
\max K_2
\approx
5.04\times10^{-3},
\]

\[
\max K_3
\approx
5.04\times10^{-3},
\]

and

\[
\max K_4
\approx
5.04\times10^{-3}.
\]

Their absolute separation is small because the shared macroscopic gradient background remains comparable to or larger than the order-specific TNE terms.

The relaxation factors nevertheless remain clearly ordered. The third- and fourth-order channels exhibit stronger kinetic activation than the second-order channel, in agreement with the prescribed relaxation spectrum.

The minimum discrete population remains positive,

\[
\min_i f_i
=
5.55\times10^{-7}.
\]

\begin{figure}[tbp]
\centering
\includegraphics[width=\textwidth,height=0.68\textheight,keepaspectratio]{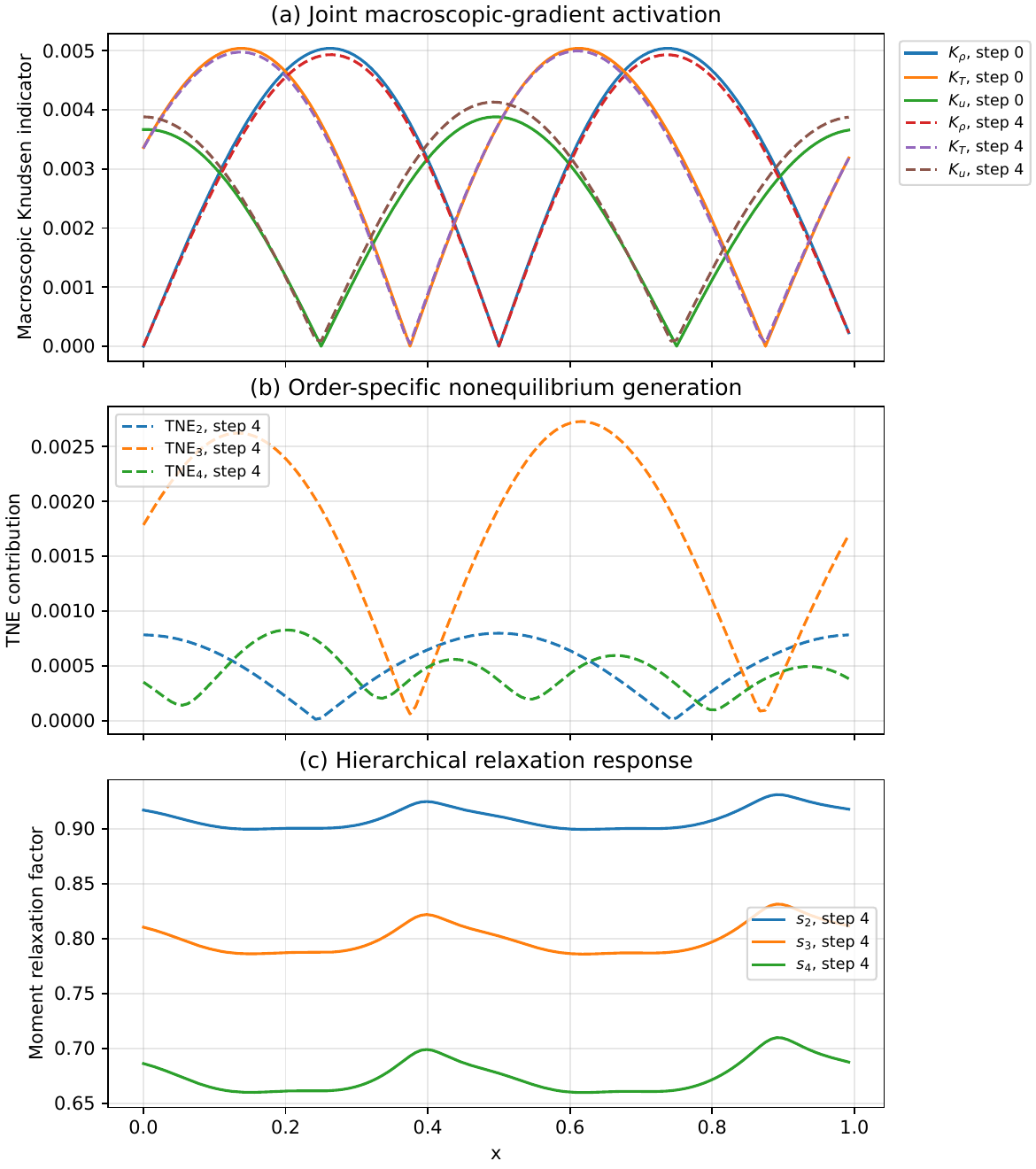}
\caption{Joint macroscopic-gradient activation and hierarchical relaxation in a smooth compression wave. (a) Density-, temperature-, and velocity-gradient rarefaction indicators at the initial state and after four timesteps. Their different spatial phases demonstrate simultaneous but noncoincident macroscopic activation. (b) Order-specific thermodynamic nonequilibrium contributions after four timesteps. The third-order contribution dominates, while the second- and fourth-order contributions remain finite and spatially structured. (c) Corresponding second-, third-, and fourth-order moment-relaxation factors. The shared macroscopic-gradient background is converted into distinct local collision responses through the order-resolved nonequilibrium indicators.}
\label{fig:smooth-compression-wave}
\end{figure}

\subsubsection{Direct common-versus-resolved comparison}

A direct comparison is first performed at

\[
\lambda=0.01.
\]

After four timesteps, the peak total TNE is

\[
3.820974\times10^{-3}
\]

for the common-sensor model and

\[
3.815802\times10^{-3}
\]

for the order-resolved model. The relative difference is therefore less than one percent.

This close agreement is expected because the shared macroscopic-gradient background dominates the effective rarefaction measure at this macroscopic-gradient sensor scale. The result is important because it shows that order-resolved activation does not introduce an artificial discrepancy in a smooth regime controlled primarily by common macroscopic rarefaction.

To expose the regime in which the hierarchical correction is strongest, the comparison is repeated in the TNE-only sensor limit,

\[
\lambda=0.
\]

In this limit, the TNE channels dominate the effective indicators. The peak total TNE decreases from

\[
3.653046\times10^{-3}
\]

in the common-sensor model to

\[
3.413224\times10^{-3}
\]

in the order-resolved model, corresponding to a reduction of

\[
6.565\%.
\]

The total-TNE difference remains negative at every spatial cell and lies in the range

\[
-2.576247\times10^{-4}
\le
\Delta\mathcal{E}_{\mathrm{tot}}
\le
-3.156659\times10^{-5},
\]

where

\[
\Delta\mathcal{E}_{\mathrm{tot}}
=
\mathcal{E}_{\mathrm{tot}}^{\mathrm{resolved}}
-
\mathcal{E}_{\mathrm{tot}}^{\mathrm{common}}.
\]

The same domain-wide reduction occurs independently in the retained moment sectors. Their difference ranges are

\[
-4.160023\times10^{-5}
\le
\Delta\mathcal{E}_2
\le
-5.069711\times10^{-7},
\]

\[
-1.025443\times10^{-4}
\le
\Delta\mathcal{E}_3
\le
-3.520555\times10^{-6},
\]

and

\[
-1.335834\times10^{-4}
\le
\Delta\mathcal{E}_4
\le
-1.372267\times10^{-5}.
\]

All three differences are negative over the complete domain. The reduction in total nonequilibrium is therefore not caused by a transfer of nonequilibrium between moment orders; all retained sectors are suppressed simultaneously.

Define the local relaxation-factor difference by

\[
\Delta s_n(x;\lambda)
=
s_n^{\mathrm{resolved}}(x;\lambda)
-
s_n^{\mathrm{common}}(x;\lambda).
\]

In the TNE-only sensor limit, this correction is positive throughout the domain for all three retained orders. The second- and third-order differences reach

\[
\max_x\Delta s_2
=
6.8944\times10^{-2}
\]

and

\[
\max_x\Delta s_3
=
6.9441\times10^{-2},
\]

respectively. The fourth-order correction is substantially larger and spans

\[
7.55\times10^{-2}
\le
\Delta s_4(x;0)
\le
2.3691\times10^{-1}.
\]

This result identifies the highest retained moment sector as the most sensitive component of the hierarchical activation mechanism.

\begin{figure}[tbp]
\centering
\includegraphics[width=\textwidth,height=0.68\textheight,keepaspectratio]{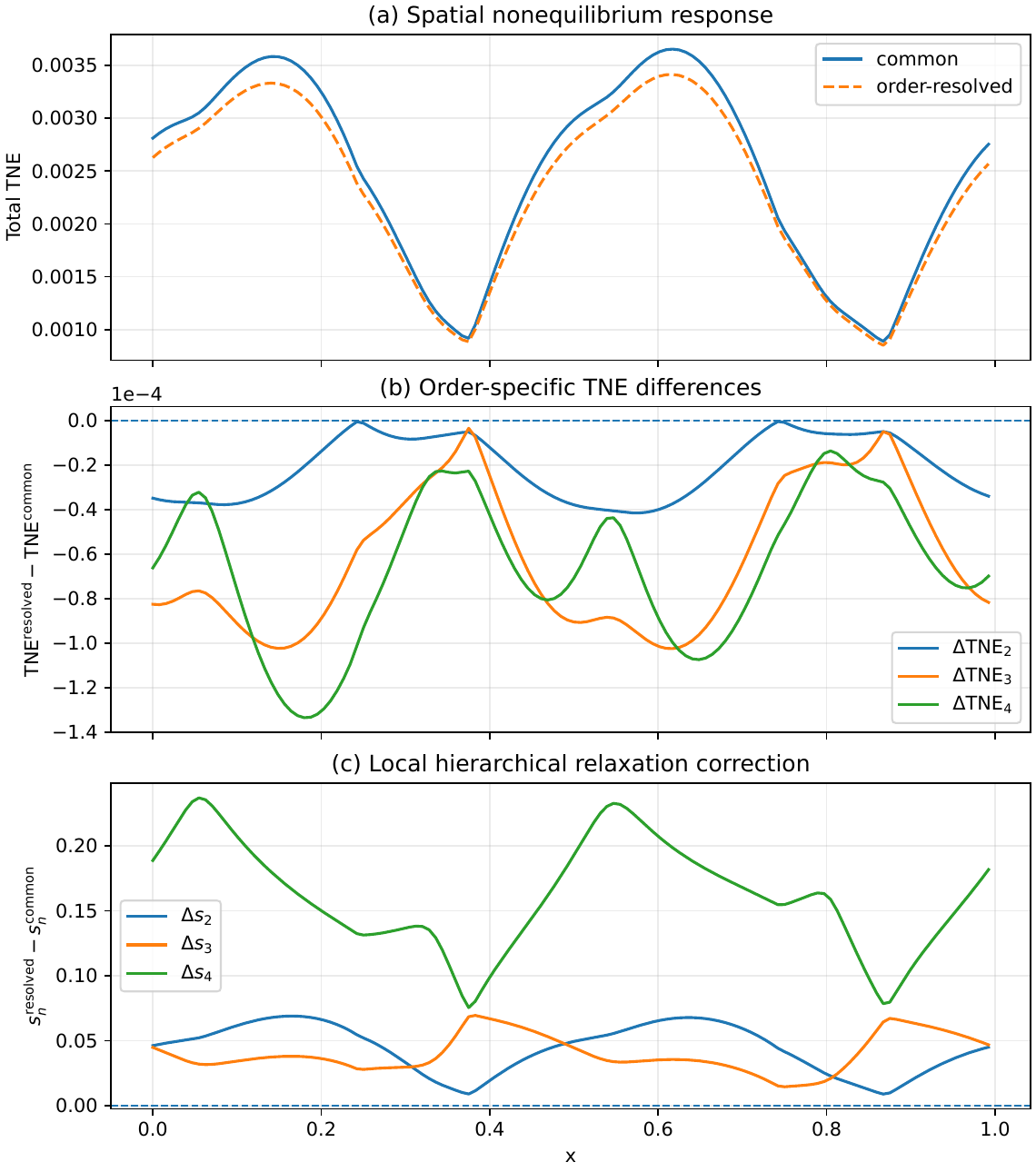}
\caption{Spatial comparison of the common and order-resolved activation models in the TNE-only sensor limit, \(\lambda=0\). (a) Total thermodynamic nonequilibrium intensity after four timesteps. The order-resolved model reduces the nonequilibrium response throughout the domain. (b) Order-specific differences in the second-, third-, and fourth-order nonequilibrium contributions. All three differences remain negative. (c) Local relaxation-factor differences \(\Delta s_n(x)=s_n^{\mathrm{resolved}}(x)-s_n^{\mathrm{common}}(x)\). The fourth-order channel exhibits the largest correction.}
\label{fig:compression-wave-model-comparison}
\end{figure}

\subsubsection{Macroscopic-gradient sensor-scale dependence}

The regime dependence is quantified by varying the macroscopic-gradient sensor scale over

\[
0\le\lambda\le0.04.
\]

The ratio of the order-resolved to common peak total TNE is

\[
R_{\mathrm{peak}}
=
\frac{
\max_x\mathcal{E}_{\mathrm{tot}}^{\mathrm{resolved}}
}{
\max_x\mathcal{E}_{\mathrm{tot}}^{\mathrm{common}}
}.
\]

In the TNE-only sensor limit,

\[
R_{\mathrm{peak}}(0)=0.934350,
\]

which corresponds to the \(6.565\%\) reduction reported above. The ratio then increases monotonically with the macroscopic-gradient sensor scale and reaches

\[
0.998646
\]

at \(\lambda=0.01\),

\[
0.999986
\]

at \(\lambda=0.02\), and

\[
0.9999999
\]

at \(\lambda=0.04\).

To summarize the local relaxation-factor differences, define

\[
D_n(\lambda)
=
\max_x
\left|
\Delta s_n(x;\lambda)
\right|.
\]

In the TNE-only sensor limit,

\[
D_2(0)=6.894368\times10^{-2},
\]

\[
D_3(0)=6.944095\times10^{-2},
\]

and

\[
D_4(0)=2.369132\times10^{-1}.
\]

At \(\lambda=0.04\), the corresponding values have decreased to

\[
D_2(0.04)=1.890487\times10^{-7},
\]

\[
D_3(0.04)=1.729505\times10^{-7},
\]

and

\[
D_4(0.04)=1.651507\times10^{-7}.
\]

The order-resolved and common formulations therefore become effectively indistinguishable when the shared macroscopic-gradient background dominates the sensor.

\begin{figure}[tbp]
\centering
\includegraphics[width=\textwidth,height=0.68\textheight,keepaspectratio]{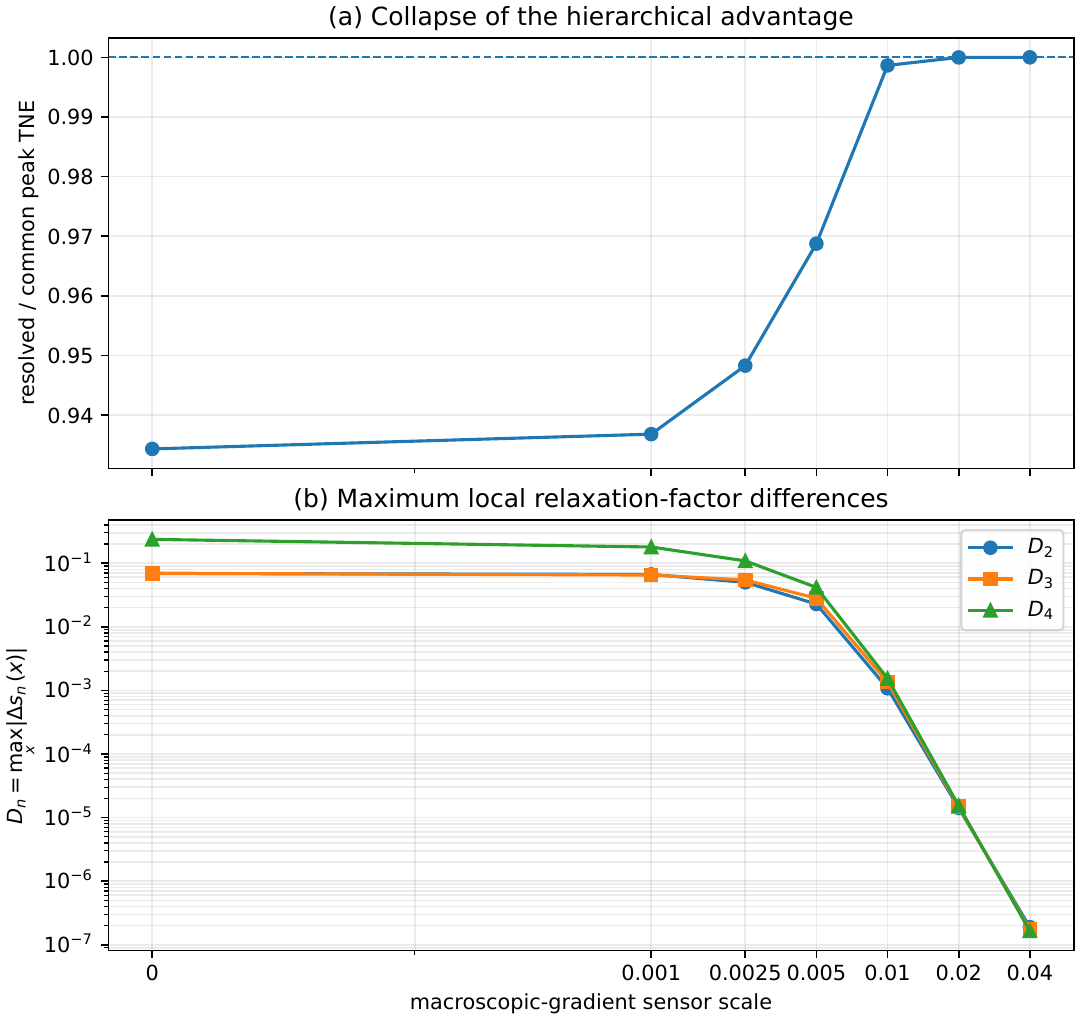}
\caption{Macroscopic-gradient sensor-scale dependence of the common and order-resolved activation models. (a) Ratio of the peak total thermodynamic nonequilibrium intensity predicted by the order-resolved model to that predicted by the common-sensor model. The ratio approaches unity as the macroscopic-gradient sensor scale increases. (b) Maximum absolute local relaxation-factor differences \(D_n(\lambda)=\max_x|\Delta s_n(x;\lambda)|\) for the second-, third-, and fourth-order sectors. These differences collapse by several orders of magnitude as the shared macroscopic-gradient background becomes dominant.}
\label{fig:compression-wave-mfp-scan}
\end{figure}

\subsubsection{Grid-and-timestep sensitivity}

A fixed-final-time grid-and-timestep sensitivity scan is performed with

\[
(N_x,N_t)
=
(32,1),
(64,2),
(128,4),
(256,8).
\]

Because the timestep scales with the cell width, all calculations reach

\[
t_f
=
4.375265\times10^{-3}.
\]

The reference-step relaxation factors are mapped to the actual timestep through exponential survival scaling. The scan therefore no longer repeats an unscaled reference-step collision factor as the timestep is refined.

The measured peak total TNE values are

\begin{table}[htbp]
\centering
\scriptsize
\setlength{\tabcolsep}{3pt}
\resizebox{\textwidth}{!}{%
\begin{tabular}{rrrrrr}
\toprule
\(N_x\) & \(N_t\) & Common model & Order-resolved model & Ratio & Reduction \\
\midrule
32 & 1 & \(1.226703\times10^{-2}\) & \(1.226703\times10^{-2}\) & 1.000000 & \(0.000\%\) \\
64 & 2 & \(6.394078\times10^{-3}\) & \(6.239484\times10^{-3}\) & 0.975822 & \(2.418\%\) \\
128 & 4 & \(3.653046\times10^{-3}\) & \(3.413224\times10^{-3}\) & 0.934350 & \(6.565\%\) \\
256 & 8 & \(2.384581\times10^{-3}\) & \(2.109795\times10^{-3}\) & 0.884766 & \(11.523\%\) \\
\bottomrule
\end{tabular}%
}
\end{table}

The one-step 32-cell calculation does not yet apply a TNE-sensitive collision to nonequilibrium generated by transport. Its initial distribution is locally equilibrated, the first collision therefore sees negligible TNE, and transport-generated nonequilibrium appears only after that collision has been evaluated.

For the three finer discretizations, the order-resolved model gives successively lower peak total TNE over the tested sequence. This trend is not interpreted as a universal asymptotic percentage or as pure spatial convergence. Spatial resolution, timestep, collision--transport splitting error, and dynamically evolving sensor fields still vary together.

\begin{figure}[tbp]
\centering
\includegraphics[width=\textwidth,height=0.68\textheight,keepaspectratio]{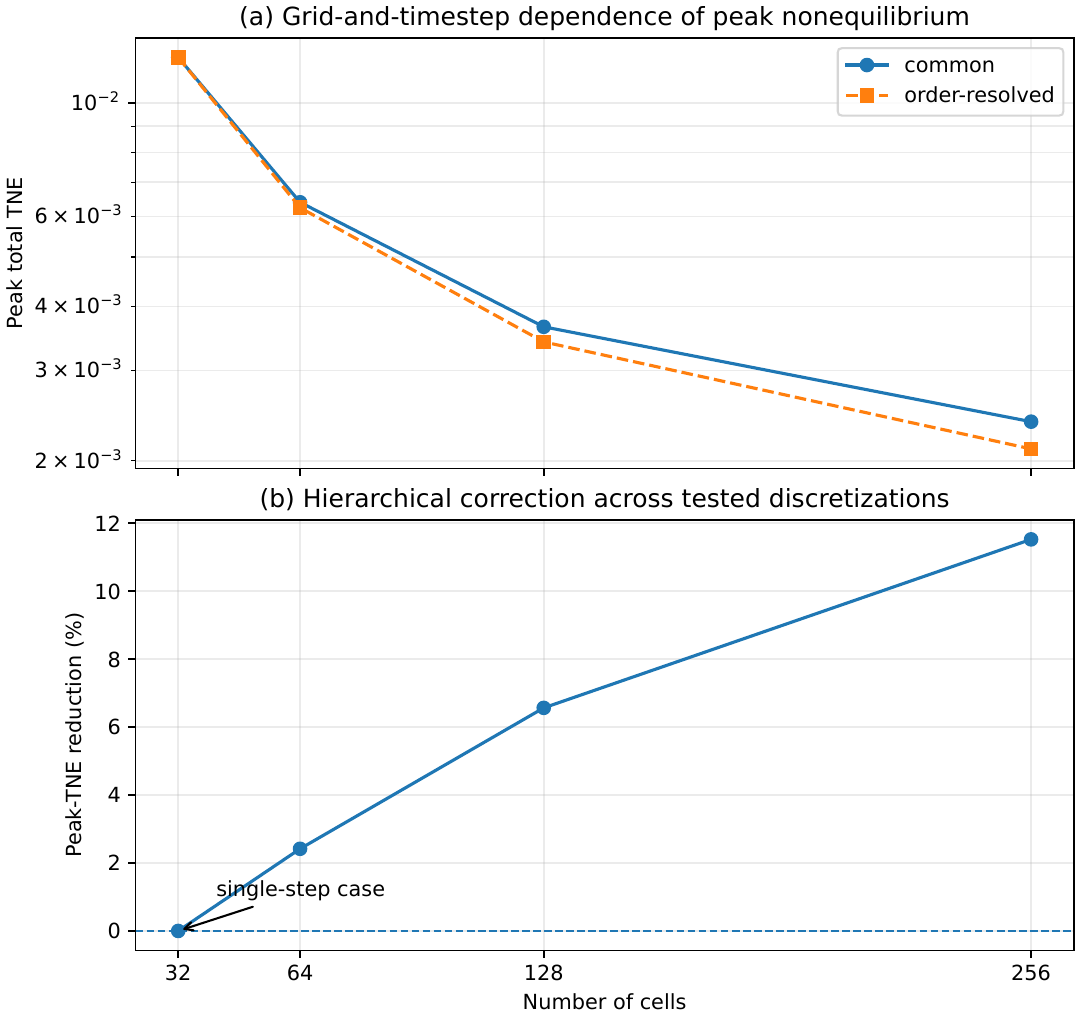}
\caption{Fixed-final-time grid-and-timestep sensitivity with timestep-consistent collision scaling. (a) Peak total thermodynamic nonequilibrium intensity predicted by the common and order-resolved models for \((N_x,N_t)=(32,1),(64,2),(128,4),(256,8)\), all at final time \(4.375265\times10^{-3}\). (b) The corresponding peak-TNE reductions are \(0.000\%\), \(2.418\%\), \(6.565\%\), and \(11.523\%\). The one-step 32-cell calculation has not yet applied a TNE-sensitive collision to transport-generated nonequilibrium. The scan retains coupled spatial, temporal, operator-splitting, and sensor-evolution effects and is not a pure spatial-convergence study.}
\label{fig:compression-wave-resolution-scan}
\end{figure}

\subsubsection{Interpretation}

The compression-wave benchmark reveals the two limiting behaviors of the hierarchical model.

When the shared macroscopic-gradient background is strong, the order-resolved indicators remain close and the model smoothly recovers the common-sensor response. When the background is weak, the distinct second-, third-, and fourth-order TNE contributions determine different local relaxation factors and produce a measurable reduction in the nonequilibrium response.

The largest correction occurs in the fourth-order channel. This is consistent with its smaller transition scale and slower kinetic-limit relaxation in the prescribed spectrum. The effect is spatially distributed and suppresses all retained nonequilibrium sectors rather than redistributing intensity among them.

The benchmark therefore provides direct evidence that hierarchical activation matters primarily in the TNE-dominated sensor regime while retaining compatibility with the common formulation in the macroscopic-gradient-dominated limit.

\subsection{Conservation diagnostics}

\begin{table}[htbp]
\centering
\scriptsize
\setlength{\tabcolsep}{3pt}
\caption{Conservation errors in the principal periodic benchmarks}
\resizebox{\textwidth}{!}{%
\begin{tabular}{llrrrr}
\toprule
Benchmark & Model & Relative mass error & Absolute momentum change & Relative total discrete kinetic-energy error & Minimum population \\
\midrule
Smooth temperature wave, \(\lambda=0.01\) & Order resolved & \num{1.221e-15} & \num{1.311e-15} & \num{7.401e-16} & \num{5.306683e-07} \\
Smooth compression wave, \(\lambda=0.01\) & Common & \num{8.882e-16} & \num{3.570e-15} & \num{7.380e-16} & \num{5.554026e-07} \\
Smooth compression wave, \(\lambda=0.01\) & Order resolved & \num{1.110e-15} & \num{3.428e-14} & \num{5.904e-16} & \num{5.554026e-07} \\
Smooth compression wave, \(\lambda=0\) & Common & \num{8.882e-16} & \num{1.697e-14} & \num{7.380e-16} & \num{5.553863e-07} \\
Smooth compression wave, \(\lambda=0\) & Order resolved & \num{9.992e-16} & \num{5.609e-15} & \num{7.380e-16} & \num{5.553837e-07} \\
\bottomrule
\end{tabular}%
}
\end{table}

For mass, the reported relative error is

\[
\varepsilon_M
=
\frac{
\left|M(t_f)-M(0)\right|
}{
\max\left(|M(0)|,10^{-30}\right)
}.
\]

The initial total momentum of the smooth periodic waves is zero up to floating-point quadrature error. A relative momentum error would therefore be poorly conditioned. The table instead reports the absolute Euclidean momentum change,

\[
\Delta P
=
\left\|
\boldsymbol P(t_f)-\boldsymbol P(0)
\right\|_2.
\]

The total discrete kinetic energy is

\[
E_{\mathrm{kin}}
=
\frac{1}{2}
\sum_{j,i}
f_{j,i}
\left\|
\boldsymbol\xi_i
\right\|_2^2,
\]

and its relative error is

\[
\varepsilon_E
=
\frac{
\left|
E_{\mathrm{kin}}(t_f)
-
E_{\mathrm{kin}}(0)
\right|
}{
\max\left(
|E_{\mathrm{kin}}(0)|,
10^{-30}
\right)
}.
\]

Because the benchmarks are periodic, the conservative finite-volume flux differences telescope globally for each discrete population. The global transport contributions to mass, momentum, and total discrete kinetic energy therefore vanish up to floating-point accumulation error; spatial discretization affects solution accuracy but not this conservation identity.

The values are computed directly from the initial and final distributions using the invariant functions implemented in the solver.

\subsection{Long-time invariant and entropy diagnostics}

A 4096-step periodic smooth-compression-wave calculation was performed on 128 cells at CFL \(=0.05\) in the TNE-only sensor limit, \(\lambda=0\). Both the common-sensor and order-resolved formulations were advanced with the same first-order upwind transport discretization.

\begin{table}[htbp]
\centering
\scriptsize
\setlength{\tabcolsep}{3pt}
\resizebox{\textwidth}{!}{%
\begin{tabular}{lrrrr}
\toprule
Model & Minimum population & Max. relative mass error & Max. absolute momentum change & Max. relative energy error \\
\midrule
Common sensor & \num{5.537336e-07} & \num{1.014e-12} & \num{2.869e-13} & \num{6.301e-13} \\
Order resolved & \num{5.537336e-07} & \num{1.012e-12} & \num{2.432e-13} & \num{6.290e-13} \\
\bottomrule
\end{tabular}%
}
\end{table}

The entropy diagnostics from the same calculations are summarized separately to avoid compressing invariant and entropy quantities into one wide table.

\begin{table}[htbp]
\centering
\small
\begin{tabularx}{\textwidth}{>{\raggedright\arraybackslash}Xrr}
\toprule
Model & Final \(H-H_0\) & Maximum absolute \(r_H\) \\
\midrule
Common & \num{-1.310743e-01} & \num{0.000e+00} \\
Order resolved & \num{-1.295946e-01} & \num{0.000e+00} \\
\bottomrule
\end{tabularx}
\end{table}

Both formulations complete the full run without nonfinite values or negative populations. The principal invariant errors remain near floating-point accumulation levels over the 4096-step sequence.

The discrete Boltzmann functional is

\[
H[f]
=
\sum_{j,i}
f_{j,i}
\ln\left(
\frac{f_{j,i}}{w_i}
\right).
\]

Its stepwise change is decomposed into collision and transport contributions,

\[
H(t)-H(0)
=
\sum_k \Delta H_{\mathrm{coll},k}
+
\sum_k \Delta H_{\mathrm{trans},k}
+
r_H,
\]

where \(r_H\) is the stagewise entropy-accounting residual. Because the collision and transport contributions are evaluated from consecutive states of the same split update, this residual primarily checks numerical bookkeeping closure. It remains at numerical zero for both models.

Neither the recorded collision-stage nor transport-stage change is positive at any completed step, and the total \(H\) decreases over the reported run. These observed signs are numerical results for this calculation only; neither the vanishing accounting residual nor the decrease of \(H\) is interpreted as a proof of a general discrete \(H\)-theorem for the combined truncated-Hermite collision, adaptive sensor, and finite-volume transport update.

\subsection{Computational cost}

\begin{table}[htbp]
\centering
\small
\caption{Common versus order-resolved runtime}
\begin{tabularx}{\textwidth}{>{\raggedright\arraybackslash}Xrr}
\toprule
Model & Median seconds per step & Relative cost \\
\midrule
Common sensor & \num{3.255574e-03} & 1.000 \\
Order resolved & \num{3.324875e-03} & 1.021 \\
\bottomrule
\end{tabularx}
\end{table}

The benchmark uses 128 cells and 20 timesteps per timing sample. Each model is measured 10 times after warmup, and the reported value is the median. The execution order is alternated between samples to reduce systematic bias from runtime drift.

The measured order-resolved runtime difference relative to the common-sensor implementation is

\[
2.13\%.
\]

These timings characterize the present NumPy implementation and should not be interpreted as hardware-independent algorithmic constants.

\subsection{Model and numerical parameters}

\begin{table}[htbp]
\centering
\scriptsize
\setlength{\tabcolsep}{3pt}
\caption{Hierarchical relaxation-spectrum parameters}
\resizebox{\textwidth}{!}{%
\begin{tabular}{rlrrrr}
\toprule
Moment order \(n\) & Primary interpretation & \(K_{0,n}\) & \(\sigma_n\) & \(s_{n,\mathrm{cont}}\) & \(s_{n,\mathrm{kin}}\) \\
\midrule
2 & Stress-like nonequilibrium & 0.050 & 2.0 & 1.00 & 0.20 \\
3 & Heat-flux-like nonequilibrium & 0.030 & 2.5 & 1.00 & 0.10 \\
4 & Higher-order kinetic nonequilibrium & 0.015 & 3.0 & 1.00 & 0.05 \\
\bottomrule
\end{tabular}%
}
\end{table}

The continuum weight is

\[
w_{c,n}(K)
=
\frac{1}{2}
\operatorname{erfc}
\left[
\frac{
\ln\!\left(
\max(K,10^{-14})/K_{0,n}
\right)
}{
\sqrt{2}\sigma_n
}
\right],
\]

and the calibrated reference-step relaxation factor is

\[
s_{n,\mathrm{ref}}
=
w_{c,n}s_{n,\mathrm{cont}}
+
(1-w_{c,n})s_{n,\mathrm{kin}}.
\]

The factor is applied directly as \(\boldsymbol{a}_{\mathrm{neq}}^{(n),*} =(1-s_n)\boldsymbol{a}_{\mathrm{neq}}^{(n)}\) at the reference timestep \(\Delta t_{\mathrm{ref}}=1.0938161705673147\times10^{-3}\). At any other positive timestep it is mapped according to \(s_n(K_n,\Delta t)=1-[1-s_{n,\mathrm{ref}}(K_n)]^{\Delta t/\Delta t_{\mathrm{ref}}}\).

\begin{table}[htbp]
\centering
\small
\caption{Default sensor and positivity parameters}
\begin{tabularx}{\textwidth}{>{\raggedright\arraybackslash}Xr>{\raggedright\arraybackslash}X}
\toprule
Parameter & Default value & Role \\
\midrule
Rarefaction norm order \(p\) & 8 & Combines gradient and TNE indicators \\
Sensor floor \(\varepsilon\) & \(10^{-14}\) & Protects normalized scales \\
Log-Gaussian floor & \(10^{-14}\) & Prevents evaluation of \(\ln 0\) \\
Positivity tolerance & \(10^{-14}\) & Separates roundoff from genuine negativity \\
Population floor & 0 & Requires nonnegative populations \\
Velocity reference scale & \(\sqrt{\max(T,10^{-14})}\) & Normalizes \(K_u\) \\
\bottomrule
\end{tabularx}
\end{table}

\begin{table}[htbp]
\centering
\small
\caption{Representation and relaxation conventions}
\begin{tabularx}{\textwidth}{>{\raggedright\arraybackslash}X>{\raggedright\arraybackslash}X>{\raggedright\arraybackslash}X}
\toprule
Parameter & Value & Interpretation \\
\midrule
Velocity set & D3Q125 & Fixed tensor-product representation \\
Highest Hermite order & 4 & Retains orders 2--4 \\
\(c_2,c_3,c_4\) & 1, 1, 1 & TNE norm multipliers \\
Tensor norm & Full Cartesian Frobenius norm & Includes permutation multiplicities \\
Sensor frame & Laboratory-frame raw Hermite sectors & Not claimed Galilean invariant \\
\(s_{n,\mathrm{ref}}\) & Reference-step discrete factor & Mapped exponentially to the actual \(\Delta t\) \\
\bottomrule
\end{tabularx}
\end{table}

\begin{table}[htbp]
\centering
\scriptsize
\setlength{\tabcolsep}{3pt}
\caption{Principal smooth-wave configurations}
\resizebox{\textwidth}{!}{%
\begin{tabular}{lrrrl}
\toprule
Benchmark & Cells & Steps & \(\lambda\) & Boundary \\
\midrule
Smooth temperature wave & 128 & 0, 1, 2, 4 & 0.01 & Periodic \\
Smooth compression wave & 128 & 0, 1, 2, 4 & 0.01 unless varied & Periodic \\
Grid-and-timestep sensitivity & 32, 64, 128, 256 & 1, 2, 4, 8 & 0 & Periodic \\
\bottomrule
\end{tabular}%
}
\end{table}

The temperature-wave initial state is

\[
\rho=1,
\qquad
u_x=0,
\qquad
T=1+0.15\cos(2\pi x/L).
\]

The compression-wave initial state is

\[
\rho=1+0.08\cos(2\pi x/L),
\]

\[
u_x=0.06\sin(2\pi x/L),
\]

and

\[
T=1+0.08\cos(2\pi x/L+\pi/4).
\]

The tabulated relaxation parameters are research-baseline values rather than universally calibrated physical constants.

\subsection{Quantitative assessment of frame dependence}

The order-resolved nonequilibrium indicators are constructed from raw Hermite coefficients defined in the fixed laboratory-frame discrete velocity basis. A uniform velocity translation therefore mixes retained orders even when the underlying physical perturbations are related by a boost. With \(\boldsymbol{U}\) denoting the boost and \(\operatorname{Sym}\) denoting normalized full index symmetrization, the deviations through fourth order obey

\[
\Delta\boldsymbol{a}^{(2)\prime}
=
\Delta\boldsymbol{a}^{(2)},
\]

\[
\Delta\boldsymbol{a}^{(3)\prime}
=
\Delta\boldsymbol{a}^{(3)}
+
3\operatorname{Sym}
\left(
\boldsymbol{U}\otimes
\Delta\boldsymbol{a}^{(2)}
\right),
\]

and

\[
\Delta\boldsymbol{a}^{(4)\prime}
=
\Delta\boldsymbol{a}^{(4)}
+
4\operatorname{Sym}
\left(
\boldsymbol{U}\otimes
\Delta\boldsymbol{a}^{(3)}
\right)
+
6\operatorname{Sym}
\left(
\boldsymbol{U}\otimes\boldsymbol{U}\otimes
\Delta\boldsymbol{a}^{(2)}
\right).
\]

Four complementary checks were used to separate this definitional dependence from evident implementation inconsistencies. First, reconstructed Maxwellians were tested over four density--temperature combinations and zero, signed, and multidirectional boosts. Their total numerical TNE remained below \(2.0\times10^{-12}\). Second, scalar local collision, loop-based batch collision, and vectorized batch collision agreed within \(2.0\times10^{-12}\) in the post-collision populations and within \(2.0\times10^{-14}\) in their reported relaxation factors. Third, second-, third-, and fourth-order perturbations over three amplitudes and three boosts reproduced the analytical translation identities within absolute tolerance \(2.0\times10^{-12}\). Fourth, homogeneous collision--projection tests recovered

\[
\Delta\boldsymbol{a}_{\mathrm{post}}^{(n)}
=
\left(1-s_n\right)
\Delta\boldsymbol{a}_{\mathrm{pre}}^{(n)}
\]

through fourth order without spatial transport.

For the baseline trace-free second-order perturbation

\[
\Delta\boldsymbol{a}^{(2)}
=
0.02
\operatorname{diag}
\left(
1,-\frac{1}{2},-\frac{1}{2}
\right),
\]

the homogeneous total TNE increased from \(0.0244949\) at \(U_0=0\) to \(0.0424524\) at \(U_0=0.2\) and \(0.0704019\) at \(U_0=0.4\). The second-order contribution and \(s_2\) remained invariant, whereas the translated third- and fourth-order contributions reduced \(s_3\) from \(1\) to \(0.570966\) and \(s_4\) from \(1\) to \(0.488813\) over the same scan. Pure third-order perturbations generated the predicted fourth-order coupling, while pure fourth-order perturbations remained invariant through the retained order.

A transport-enabled study added a constant streamwise offset \(U_0\) to the smooth compression-wave initial condition while preserving the density, temperature, and velocity-gradient profiles. Relative to \(U_0=0\), the peak total TNE increased by \(3.906\%\), \(10.474\%\), \(23.879\%\), and \(50.627\%\) for \(U_0=0.025\), \(0.05\), \(0.10\), and \(0.20\), respectively. The peak second-order contribution changed by less than \(0.1\%\), the third-order contribution changed by at most approximately \(6.1\%\), and the fourth-order increase reached approximately \(36.7\%\), \(74.5\%\), \(151\%\), and \(303\%\).

At matched final time, the relative total-TNE change for \(U_0=0.1\) was \(0.148742\) on 24 cells at CFL \(=0.4\), \(0.174160\) on 48 cells at CFL \(=0.4\), and \(0.141330\) on 24 cells at CFL \(=0.2\). Thus transport adds a measurable grid- and timestep-dependent contribution to the laboratory-frame sensitivity.

\begin{figure}[tbp]
\centering
\includegraphics[width=\textwidth,height=0.68\textheight,keepaspectratio]{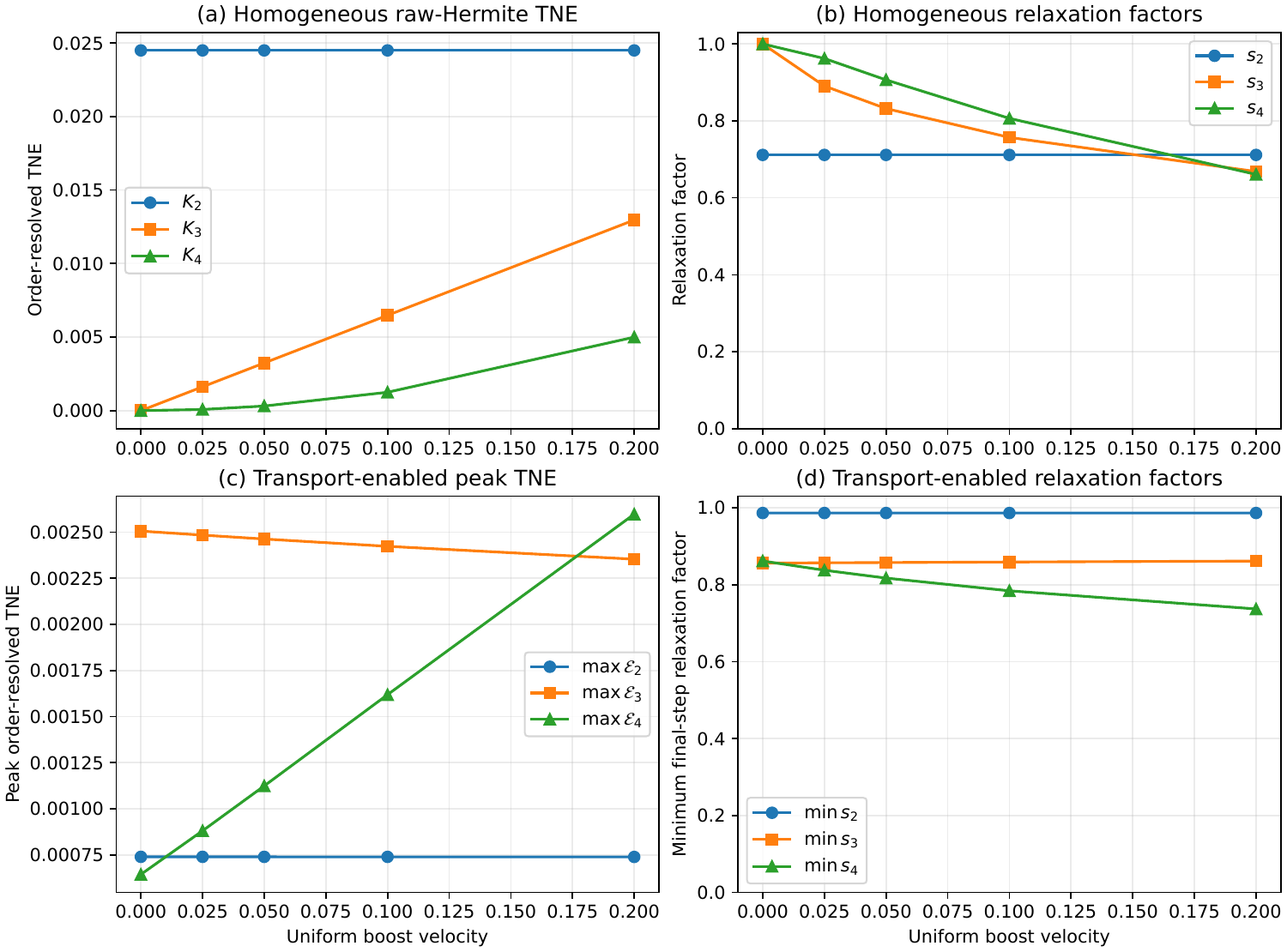}
\caption{Quantitative uniform-boost assessment of raw-Hermite TNE and relaxation. (a,b) A homogeneous trace-free second-order perturbation exhibits the analytical coupling into third- and fourth-order raw-Hermite sectors and the corresponding change in relaxation factors. (c,d) The transport-enabled smooth compression wave adds discretization-dependent changes to the peak TNE contributions and final-step relaxation factors.}
\label{fig:frame-dependence-validation}
\end{figure}

\begin{table}[htbp]
\centering
\small
\caption{Frame-dependence validation summary}
\begin{tabularx}{\textwidth}{>{\raggedright\arraybackslash}Xc>{\raggedright\arraybackslash}X}
\toprule
Validation & Result & Interpretation \\
\midrule
Boosted Maxwellian TNE & \(<2.0\times10^{-12}\) & No spurious equilibrium TNE detected \\
Scalar, loop-batch, vectorized agreement & \(<2.0\times10^{-12}\) & No evident path or broadcasting inconsistency \\
Raw-Hermite translation & \(<2.0\times10^{-12}\) absolute error & Numerical projection matches the analytical frame transformation \\
Homogeneous relaxation & \(<3.0\times10^{-12}\) absolute error & Collision applies the reported order-specific factor \\
Transport boost test & Grid/CFL dependent & Physical-space discretization adds frame sensitivity \\
\bottomrule
\end{tabularx}
\end{table}

Within the tested parameter range, no evident implementation inconsistency was detected. The dominant measured frame dependence instead follows the analytical inter-order translation structure of laboratory-frame raw Hermite moments. This statement does not exclude implementation errors outside the tested range and does not establish that a central-moment replacement would automatically remove finite-quadrature, truncation, or transport-discretization effects.

\subsection{Fixed-grid timestep sensitivity}

The fixed-final-time grid-and-timestep study shows an increasing peak-TNE reduction across the three finer tested pairs after timestep-consistent collision scaling is introduced. Because that study changes the grid, timestep, splitting error, and dynamically evaluated sensors together, it remains a coupled numerical sensitivity result rather than a pure spatial-convergence sequence.

A separate fixed-grid, fixed-final-time CFL scan was used to isolate the remaining timestep dependence at \(N_x=128\). The CFL number was reduced from \(0.4\) to \(0.0015625\), while the number of timesteps was increased from \(4\) to \(1024\) so that the final time remained unchanged. The peak total TNE decreased systematically in both formulations. Between the two finest tested timesteps, the common-sensor peak changed by \(0.654\%\), the order-resolved peak changed by \(0.819\%\), and the hierarchical peak-TNE reduction changed by \(0.126\) percentage points.

The successive timestep differences decrease approximately by a factor of two under timestep halving, consistent with first-order temporal behavior of the sequential collision--transport splitting. The hierarchical reduction increases from \(6.565\%\) at CFL \(=0.4\) to \(22.493\%\) at CFL \(=0.0015625\). Therefore, the \(6.565\%\) value characterizes the stated benchmark timestep rather than a timestep-converged universal percentage. The finest scan values provide evidence of practical timestep convergence under the adopted criteria, but they remain tied to the fixed spatial grid, transport discretization, sensor definition, and final time.

\subsection{Transport-discretization sensitivity}

A transport-discretization sensitivity study compared the baseline first-order upwind update with the available MUSCL reconstruction at fixed grid, fixed final time, and three CFL values. Such sensitivity is expected in finite-volume discrete-kinetic schemes because the physical-space reconstruction contributes its own numerical dissipation and truncation error \cite{peng1999finitevolume,patil2009tvd}. In every tested case, the order-resolved formulation retained a lower peak total TNE than the common-sensor formulation. The hierarchical reductions predicted by upwind and MUSCL were \(6.5650\%\) and \(6.5069\%\) at CFL \(=0.4\), \(18.9442\%\) and \(18.7022\%\) at CFL \(=0.05\), and \(21.6251\%\) and \(21.3404\%\) at CFL \(=0.0125\), respectively.

The absolute peak-TNE values differed by approximately \(0.5\%\) to \(1.0\%\) between the two transport schemes over the tested configurations, while the hierarchical reduction differed by less than \(0.3\) percentage points. All final distributions retained positive populations. These results show that the sign and qualitative magnitude of the hierarchical correction are not artifacts of the first-order upwind discretization, although the reported numerical values retain a finite dependence on the selected transport scheme.

\subsection{Relaxation-spectrum sensitivity}

A relaxation-spectrum sensitivity study was performed at fixed grid, fixed final time, CFL \(=0.05\), and zero macroscopic-gradient sensor scale. The transition locations \(K_{0,n}\), transition widths \(\sigma_n\), and kinetic-limit reference factors \(s_{n,\mathrm{kin}}\) were perturbed coherently around the default spectrum.

The order-resolved formulation retained a lower peak total TNE than the common-sensor formulation for every tested spectrum. The hierarchical reduction varied from \(16.4679\%\) to \(21.4968\%\), compared with \(18.9442\%\) for the baseline spectrum. Scaling the transition locations produced the largest quantitative response: multiplying all \(K_{0,n}\) by \(0.5\) increased the reduction to \(21.4968\%\), whereas multiplying them by \(2\) reduced it to \(16.4679\%\). The corresponding ranges for the transition-width and kinetic-limit-factor perturbations were narrower.

All tested final distributions retained positive populations. The sign of the hierarchical correction is therefore robust under these coordinated parameter perturbations, although its numerical magnitude remains dependent on the selected relaxation spectrum. The default spectrum should consequently be interpreted as a research baseline rather than a uniquely calibrated or universal parameter set.

\subsection{Shear-wave modal-decay and effective-viscosity refinement}

A transverse shear-wave calculation was used to test analytic modal decay, linear-response behavior, and grid refinement of the fitted effective viscosity. The initial state had uniform density and temperature and a small sinusoidal transverse velocity, \(u_y=A\sin(2\pi x/L)\). A constant effective Knudsen indicator was prescribed so that the relaxation spectrum remained fixed during the comparison.

For amplitudes \(A=0.0025\), \(0.005\), and \(0.01\), the fundamental Fourier mode exhibited near-exponential decay. The fitted decay curves had \(R^2>0.99997\) on the 64-, 128-, and 256-cell grids, with phase drift at the level of floating-point roundoff. At each resolution, the relative spread of the fitted effective viscosity over the three amplitudes was below \(0.0012\%\), confirming linear-response behavior over the tested range.

A fixed-final-time refinement using 64, 128, 256, and 512 cells gave fitted viscosities of \(6.633857\times10^{-4}\), \(5.196832\times10^{-4}\), \(4.852919\times10^{-4}\), and \(4.768346\times10^{-4}\), respectively, for \(A=0.005\). The three-level empirical convergence orders were \(2.063\) and \(2.024\). Richardson extrapolation from the finest three grids gave \(\nu_\infty=4.740766\times10^{-4}\), with the 512-cell result differing from this estimate by \(0.5818\%\).

All shear-wave calculations retained positive populations, negligible phase drift, and global invariant changes near floating-point accumulation levels. This test confirms the expected single-mode decay and provides a quantitative grid-refinement assessment of the fitted effective viscosity. It is not, however, a substitute for validation against an independent discrete-velocity, Boltzmann, or DSMC reference solution.

\section{Discussion}

\subsection{Physical interpretation of hierarchical activation}

The numerical experiments show that the principal effect of the order-resolved formulation is not the introduction of a new macroscopic rarefaction scale, but the removal of an artificial coupling among nonequilibrium moment sectors.

In the common-sensor model, the second-, third-, and fourth-order relaxation factors are all controlled by one scalar indicator. A large nonequilibrium contribution in one order therefore affects the relaxation response of the other retained orders, even when their local deviations from equilibrium are much smaller.

The hierarchical formulation replaces this shared nonequilibrium input with the three order-specific quantities

\[
\mathcal{E}_2,\qquad
\mathcal{E}_3,\qquad
\mathcal{E}_4.
\]

The macroscopic-gradient background remains common, but each retained moment sector responds to its own nonequilibrium state. This separation preserves the physical distinction among stress-like, heat-flux-like, and higher-order kinetic deviations.

The pure-order activation tests confirm that this construction is selective: a perturbation confined to one moment order activates its matching relaxation channel while leaving the nonmatching sectors at roundoff level. The mixed-order tests then show that this selectivity changes the cumulative relaxation pathway when several nonequilibrium orders coexist.

\subsection{Sensitivity of the fourth-order sector}

Across the homogeneous and spatial benchmarks, the fourth-order sector shows the largest difference between the common and order-resolved models. In the TNE-only-sensor-limit compression-wave comparison, the fourth-order factor correction reaches

\[
\Delta s_4
\approx
2.37\times10^{-1},
\]

which is substantially larger than the maximum second- and third-order corrections.

This stronger response follows from the prescribed relaxation spectrum for each moment order. The fourth-order channel has the smallest transition scale and the slowest kinetic-limit relaxation factor. It therefore reacts more strongly to a separation between the common total-TNE indicator and the local fourth-order TNE contribution.

This result does not imply that fourth-order nonequilibrium always dominates the total kinetic response. In the temperature- and compression-wave benchmarks, the third-order TNE signal is frequently the largest contribution. Rather, the fourth-order sector is the most sensitive to the choice of activation model because its relaxation curve is intentionally more kinetic.

\subsection{Regime dependence}

The macroscopic-gradient sensor-scale scan identifies two distinct regimes.

For small \(\lambda\), the shared gradient contribution is weak and the effective indicators are controlled mainly by the order-specific TNE terms. Differences among \(\mathcal{E}_2\), \(\mathcal{E}_3\), and \(\mathcal{E}_4\) then produce distinct local relaxation factors and a measurable change in the nonequilibrium response.

At zero macroscopic-gradient sensor scale, corresponding to the TNE-only sensor limit, the order-resolved model lowers the peak total TNE in the smooth compression wave by

\[
6.565\%.
\]

The reduction occurs at every spatial cell and in every retained moment sector.

As \(\lambda\) increases, the shared gradient background becomes dominant. The three effective indicators then satisfy approximately

\[
K_2\approx K_3\approx K_4,
\]

and the common and order-resolved formulations converge.

At

\[
\lambda=0.04,
\]

the peak-TNE ratio is indistinguishable from unity at the plotted precision, and the three maximum absolute local relaxation-factor differences \(D_n\) are of order

\[
10^{-7}.
\]

This limiting behavior shows that hierarchical activation does not force a persistent discrepancy between the two models. The order-resolved formulation changes the dynamics only when the nonequilibrium channels contain information that is lost by scalar aggregation.

\subsection{Interpretation of the reduced TNE response}

The order-resolved model consistently produces a smaller residual or peak total TNE in the reported comparisons. This reduction should be interpreted carefully.

The result does not imply that smaller TNE is universally more accurate or physically preferable. TNE is a diagnostic of departure from local equilibrium, not an error norm with respect to an external reference solution.

The present results demonstrate that the common-sensor formulation can retain additional nonequilibrium because a strongly activated channel reduces the relaxation of other moment sectors. The hierarchical model removes this cross-order suppression and permits weakly nonequilibrium sectors to relax more rapidly.

Whether this revised relaxation pathway improves agreement with the underlying kinetic equation must ultimately be evaluated against independent reference solutions, since high-order lattice kinetic accuracy depends on both velocity-space representation and spatial-temporal discretization \cite{shi2021accuracy}. The present study establishes the mechanism, selectivity, and numerical consequences of the model, but it does not yet claim universal accuracy superiority.

\subsection{Numerical robustness}

The implementation remains stable and positive for all reported mechanism and smooth-wave tests. The minimum populations stay strictly above zero, typically near

\[
5\times10^{-7},
\]

and no positivity limiter is required in the principal benchmarks.

The common and order-resolved models use the same discrete velocity set, equilibrium reconstruction, transport discretization, timestep, and boundary conditions. Their differences can therefore be attributed to the order-specific sensor and relaxation response rather than to a change in the underlying numerical representation.

A separate 4096-step calculation on 128 cells at CFL \(=0.05\) in the TNE-only sensor limit was completed by both formulations without negative or nonfinite populations. The maximum relative mass and energy errors remained approximately \(10^{-12}\) and \(6.3\times10^{-13}\), respectively. The stagewise entropy-accounting residual remained at numerical zero, confirming closure of the recorded split-step bookkeeping. The completed positive run and small invariant drift supply a longer-time robustness check beyond the short smooth-wave comparisons.

With timestep-consistent exponential collision scaling, the fixed-final-time grid-and-timestep scan gives reductions of \(0.000\%\), \(2.418\%\), \(6.565\%\), and \(11.523\%\) for 32, 64, 128, and 256 cells, respectively. The zero difference in the single-step 32-cell case occurs because transport-generated nonequilibrium is not encountered by a subsequent collision. Although the three finer calculations show an increasing correction over the tested sequence, the scan still changes the grid, timestep, collision--transport splitting error, and sensor evolution together. It therefore does not establish pure spatial convergence or a universal converged percentage.

\subsection{Conservation and entropy diagnostics}

The collision operator is designed to leave the conserved moments unchanged, and the solver separately records collision- and transport-stage changes in mass, momentum, and kinetic energy.

For periodic benchmarks, the conservative finite-volume flux differences telescope globally for each discrete population. The corresponding transport-stage changes in total mass, momentum, and discrete kinetic energy therefore vanish up to floating-point accumulation error. Spatial discretization affects solution accuracy but not this global conservation identity. These diagnostics provide a direct check that the adaptive relaxation mechanism does not inadvertently alter conserved quantities.

The discrete entropy-like functional is monitored only when all populations are positive. Its use in this study is diagnostic. No general discrete \(H\)-theorem is claimed for the fourth-order truncated Hermite reconstruction, order-dependent collision spectrum, finite-volume transport, and optional positivity limiter considered together.

A more complete entropy analysis would require either a collision law constructed from a discrete convex entropy, as in entropic lattice Boltzmann formulations \cite{succi2002htheorem,ansumali2002entropic}, or a separate proof of monotonicity for the full update.

The global conservation tests confirm that the collision--transport implementation preserves mass and total discrete kinetic energy to floating-point accuracy in the periodic benchmarks, while the total momentum change remains correspondingly small. All reported final states retain strictly positive populations. The detailed numerical values are collected in Table 1.

\subsection{Computational cost}

The order-resolved model reuses the same macroscopic state, Hermite coefficients, and gradient calculations required by the common-sensor formulation. Its additional work consists mainly of retaining three TNE components, constructing three effective Knudsen indicators, and evaluating three relaxation curves instead of broadcasting one common indicator.

The asymptotic cost remains proportional to the number of cells and discrete velocities. Because the velocity set is unchanged, the additional memory and runtime are expected to be modest relative to the cost of projection, reconstruction, and transport.

The timing comparison in Table 2 shows that the additional order-specific sensor and relaxation evaluations contribute only a small runtime increase in the present NumPy implementation. This measurement is implementation- and hardware-dependent, but it indicates that hierarchical activation does not dominate the total collision--transport cost.

\subsection{Limitations}

The present study has several limitations.

First, the velocity set is fixed at D3Q125. The method adapts the local collision response but does not adapt the extent or resolution of velocity space.

Second, the spatial demonstrations are one-dimensional. Although the velocity representation and moment tensors are three-dimensional, the reported transport problems do not test genuinely multidimensional gradients, vortical structures, or complex wall geometry.

Third, the benchmarks are smooth and remain positive without strong limiting. The behavior of the model near shocks, discontinuities, or strongly under-resolved kinetic layers has not yet been established.

Fourth, the present comparisons are internal comparisons between a common-sensor and an order-resolved formulation. No direct benchmark against a high-resolution discrete-velocity method, DSMC solution, or analytic kinetic reference has yet been included.

Fifth, the relaxation-spectrum parameters are prescribed rather than calibrated from first-principles transport coefficients or optimization against reference data. Their ordering is physically motivated, but their numerical values should not be interpreted as universal.

The additional sensitivity studies reported in Sections 5.10--5.14 expose three further limitations. First, the present hierarchical indicators are based on laboratory-frame raw Hermite coefficients and therefore retain a measurable dependence on a uniform velocity translation. Dedicated Maxwellian, cross-path, pure-order, and homogeneous-collision checks show that this dependence follows the analytical inter-order translation structure of raw Hermite moments, with the largest transport-enabled response in the fourth-order contribution. A locally centered or central-Hermite formulation is a natural extension, but its behavior under finite quadrature, fourth-order truncation, and physical-space transport discretization remains to be established.

Second, the fixed-grid timestep, transport-discretization, and relaxation-spectrum studies show that the sign of the hierarchical correction is robust over the tested configurations, but its numerical magnitude remains dependent on timestep, transport scheme, and prescribed spectrum. Reported percentages must therefore be interpreted as benchmark-specific rather than universal constants.

Third, the shear-wave experiment confirms the expected analytic mode and provides a refinement assessment of the fitted effective viscosity, but it does not compare that fitted value with an independently predicted transport coefficient and is not an independent kinetic-reference validation. Comparison with a separate discrete-velocity, Boltzmann, or DSMC solution remains necessary.

\subsection{Scope of future work}

The next stage of development should separate into two directions.

For the present fixed-velocity model, the highest priority is external validation. Suitable tests include comparison with a reference discrete-velocity solver, analysis of transport coefficients, and multidimensional smooth and weakly discontinuous flows.

A separate future contribution should address dynamic velocity-space adaptation. This would include switching among D3Q125, D3Q343, and D3Q729, local velocity scaling, active-node selection, and conservative remapping between velocity sets.

Keeping these two contributions separate is important. The present work isolates the effect of hierarchical moment activation, whereas dynamic velocity adaptation introduces additional approximation, conservation, and remapping questions.

\section{Conclusions}

This work introduced an order-resolved hierarchical relaxation model for a fixed D3Q125 discrete-velocity formulation. The central idea is to retain a shared macroscopic rarefaction background while allowing the second-, third-, and fourth-order moment sectors to respond to their own thermodynamic nonequilibrium indicators.

The resulting effective Knudsen indicators,

\[
K_2,\qquad K_3,\qquad K_4,
\]

drive distinct order-dependent log-Gaussian relaxation factors. This construction removes the artificial cross-order coupling introduced when all retained moments share a single scalar nonequilibrium sensor.

The pure-order activation tests verify that each perturbation activates its matching moment channel without measurable leakage into the nonmatching sectors. Homogeneous mixed-order relaxation, amplitude scans, and composition scans further show that the order-resolved model consistently produces a lower residual nonequilibrium response than the common-sensor formulation across a broad range of perturbation strengths and moment compositions.

The smooth temperature-wave benchmark demonstrates that the order-specific TNE signals are generated dynamically from an initially equilibrium distribution. The third-order contribution becomes dominant, while the prescribed relaxation spectrum maintains clearly separated second-, third-, and fourth-order collision responses.

The smooth compression-wave benchmark provides the clearest comparison between the common and order-resolved models. In the TNE-only sensor limit, \(\lambda=0\), the hierarchical formulation reduces the peak total TNE by

\[
6.565\%,
\]

and the reduction occurs at every spatial location and in every retained moment sector. The largest relaxation correction appears in the fourth-order channel.

A macroscopic-gradient sensor-scale scan shows that this distinction is regime dependent. As the shared macroscopic-gradient contribution becomes dominant, the order-resolved and common formulations converge smoothly. At \(\lambda=0.04\), the maximum absolute local relaxation-factor differences \(D_n\) are of order \(10^{-7}\), and the peak-TNE ratio is effectively unity.

With timestep-consistent exponential collision scaling, the fixed-final-time grid-and-timestep sensitivity study gives peak-TNE reductions of \(0.000\%\), \(2.418\%\), \(6.565\%\), and \(11.523\%\) for 32, 64, 128, and 256 cells. The single-step 32-cell calculation does not expose transport-generated nonequilibrium to a subsequent TNE-sensitive collision. The scan still couples spatial resolution, timestep, operator splitting, and sensor evolution, so it is not interpreted as a pure spatial-convergence test or as evidence for a universal converged percentage.

All principal benchmarks remain positive without requiring the positivity limiter. In a separate 4096-step calculation on 128 cells, both formulations remain positive and finite, with maximum relative mass and energy errors of approximately \(10^{-12}\) and \(6.3\times10^{-13}\), respectively, while the stagewise entropy-accounting residual remains at numerical zero. This residual verifies closure of the recorded split-step bookkeeping and is not used as evidence of a general discrete \(H\)-theorem.

The present results establish the selectivity, regime dependence, and numerical consequences of order-resolved moment activation on a fixed D3Q125 velocity set. They do not yet establish universal accuracy superiority over the common-sensor formulation. The next priority is external validation against independent kinetic reference solutions and extension to multidimensional flows.

The supplementary sensitivity studies show that the hierarchical correction retains its sign under timestep refinement, transport-scheme changes, and coordinated relaxation-spectrum perturbations, although its numerical magnitude remains configuration-dependent. The uniform-boost study also identifies the laboratory-frame raw-Hermite sensor as a source of frame dependence, especially in the fourth-order sector.

Quantitative uniform-boost tests further identify the measurable frame sensitivity of the raw-Hermite indicators and separate it from evident collision-implementation inconsistencies within the tested parameter range.

The shear-wave calculation confirms near-exponential modal decay, linear-response behavior, and approximately second-order refinement of the fitted effective viscosity. Independent kinetic-reference validation and genuinely multidimensional tests nevertheless remain necessary. Dynamic velocity-space adaptation is a separate future direction.

The periodic benchmarks preserve the principal global invariants to floating-point accuracy and retain positive populations. In the reference NumPy implementation, the order-resolved formulation adds about \(2.13\%\) to the measured runtime in the reported compression-wave timing test.

\bibliographystyle{plain}
\bibliography{references}

\end{document}